\documentclass[12pt,a4paper]{article}
\usepackage{latexsym}
\usepackage{mathrsfs}
\usepackage{amsmath}
\usepackage{amssymb}
\usepackage[usenames,dvipsnames,svgnames,table]{xcolor}

\newtheorem{theorem}{Theorem}[section]
\newtheorem{lemma}[theorem]{Lemma}

\def \D {{\mathcal{D}}}

\def \P {{\mathcal {P}}}

\def \E {\varepsilon}

\def \l {\lambda}
\def \e {\E}

\def \deg {{\rm deg}}

\def \proof {\noindent{\bf Proof}\quad}
\def \binom#1#2{{#1\choose#2}}

\newcommand{\qed}{\hfill$\Box$\vspace{0.2cm}}

\title{\bf Decompositions of complete multigraphs into cycles of varying lengths}

\author{
Darryn Bryant\thanks{School of Mathematics and Physics, The University of Queensland, QLD 4072, Australia.\newline Email: \texttt{db@maths.uq.edu.au}, \qquad \texttt{bmm@maths.uq.edu.au},\qquad  \texttt{bsmith.maths@gmail.com}}, Daniel Horsley\thanks{School of Mathematical Sciences, Monash University, VIC 3800, Australia. \newline  Email: \texttt{daniel.horsley@monash.edu}}, Barbara Maenhaut\footnotemark[1]\, and \\Benjamin R. Smith\footnotemark[1]
}

\date{ }

\begin{document}
\maketitle\thispagestyle{empty}
\def\baselinestretch{1}\small\normalsize

\begin{abstract}
We establish necessary and sufficient conditions for the existence of 
a decomposition of a complete multigraph into edge-disjoint cycles of specified lengths, 
or into edge-disjoint cycles of specified lengths and a perfect matching. 
\end{abstract}

\section{Introduction}

A \emph{decomposition} of a graph $K$ is a collection $\D$ of subgraphs of $K$ such that the edge sets of the graphs in $\D$ partition the edge set of $K$.
If the complete graph of order $n$, denoted $K_n$, has a decomposition into $t$ 
cycles of specified lengths $m_1,m_2,\ldots,m_t$, then it is easy to see that 
$3\leq m_i\leq n$ for $i=1,2,\ldots,t$, $n$ is odd, and $m_1+m_2+\cdots+m_t=\binom n2$.
Similarly, if $K_n$ has a decomposition into $t$
cycles of specified lengths $m_1,m_2,\ldots,m_t$ and a perfect matching, then 
$3\leq m_i\leq n$ for $i=1,2,\ldots,t$, $n$ is even, and $m_1+m_2+\cdots+m_t=\binom n2-\frac n2$.
In \cite{BrHoPe} it was shown that these obvious necessary conditions are also sufficient for 
the existence of the desired decomposition, thereby 
solving a problem posed by Alspach in 1981 \cite{Als}.

In this paper, the analogous problem for decompositions of complete multigraphs into 
cycles of specified lengths is completely solved, see Theorem \ref{mainfull}. 
The complete multigraph of order $n$ and multiplicity
$\l$, which has $\l$ distinct edges joining each pair of distinct vertices, 
is denoted $\l K_n$. 

\begin{theorem}\label{mainfull}
There is a decomposition $\{G_1,G_2,\ldots,G_t\}$ of $\lambda K_n$ in which $G_i$ is an $m_i$-cycle for $i=1,2,\ldots,t$ if and only if
\begin{itemize}
\item $\lambda (n-1)$ is even;
\item $2\leq m_1,m_2,\ldots,m_t\leq n$;
\item $m_1+m_2+\cdots+m_t=\lambda \binom{n}{2}$;
\item $\max(m_1,m_2,\ldots,m_t)+t-2\leq \frac{\lambda}{2}\binom{n}{2}$ when $\lambda$ is even; and
\item $\sum_{m_i=2} m_i\leq (\lambda-1)\binom{n}{2}$ when $\lambda$ is odd.
\end{itemize}
There is a decomposition $\{G_1,G_2,\ldots,G_t,I\}$ of $\lambda K_n$ in which $G_i$ is an $m_i$-cycle for $i=1,2,\ldots,t$ and $I$ is a perfect matching if and only if
\begin{itemize}
\item $\lambda (n-1)$ is odd;
\item $2\leq m_1,m_2,\ldots,m_t\leq n$;
\item $m_1+m_2+\cdots+m_t=\lambda \binom{n}{2}-\frac{n}{2}$; and
\item $\sum_{m_i=2} m_i\leq (\lambda-1)\binom{n}{2}$.
\end{itemize}
\end{theorem}

The necessity of the conditions of Theorem \ref{mainfull} are proved in Section \ref{Section2} and sufficiency is 
proved in Section \ref{SectionProof}.
Note that for $\l=1$, the condition that $\sum_{m_i=2} m_i\leq (\lambda-1)\binom{n}{2}$ implies that each 
$m_i\geq 3$, and so the necessary conditions of Theorem \ref{mainfull} reduce to the familiar necessary 
conditions for Alspach's problem, as described in the first paragraph.  

There has been considerable work done on the existence of decompositions of complete multigraphs into
cycles. For the case $\l=1$, the eventual complete 
solution \cite{BrHoPe} was preceded by numerous partial results, dating back to 1847 \cite{Kir,Luc}.
The special case where all of the cycles have uniform length was settled by Alspach, Gavlas and \v{S}ajna
\cite{AlsGav,Saj}. Important preliminary results that contributed directly to the complete solution for $\l=1$
are given in \cite{BryHorPack,BryHorLong,BryHorAsym,BryHorMae}, and other partial results can be found in the large number of cited papers in \cite{BrHoPe} and the 
surveys \cite{Bry2} and \cite{BryRod}.

For $\l>1$, there have been relatively few results for cases where there are 
cycles of varying lengths \cite{BrHoMaSm}, but the case where all the cycles are of 
uniform length $m$ has been studied 
extensively. Solutions for small values of $m$ are given in \cite{BS,BHS,Hanani,HR,RH}, other partial results appear 
in \cite{Smith,Smith2}, and a complete solution for all $m$ and all $\l$ is given in \cite{BrHoMaSm}.

Analagous problems concerning decompositions of $\l K_n$ into paths, matchings or stars 
of sizes $m_1,m_2,\ldots,m_t$ have also been considered. The problem is completely solved for paths in 
\cite{Bry3}. Baranyai's Theorem \cite{Bar} settles the problem for decompositions of $\l K_n$ into matchings when 
$\l=1$, and an easy induction extends this result to a complete solution for all $\l\geq 1$.
The problem of decomposing $\l K_n$ into isomorphic stars has been solved \cite{TarsiStars}, as has the problem of decomposing $K_n$ into stars of sizes $m_1,m_2,\ldots,m_t$ \cite{LinShy}. 

We briefly mention some basic graph theory terminology that we will use. 
A graph $G$ is a nonempty set
$V(G)$ of vertices and a set $E(G)$ of edges, together with a function which maps each edge in
$E(G)$ to a pair of distinct vertices in $V(G)$ called its \emph{endpoints}. The {\em size} of a
graph $G$ is $|E(G)|$, and the number of edges in $G$ which have
$u$ and $v$ as their endpoints is denoted by $\mu_G(uv)$ and called the \emph{multiplicity} of edge $uv$.
Note that the above definition of a graph distinguishes different edges with the same endpoints.
Much of the time, however,
distinguishing edges with the same endpoints is an unnecessary complication which we ignore. 
For example, if we wish to delete an edge with endpoints $u$ and $v$ from some graph $G$, then
it will generally not matter which of the $\mu_G(uv)$ edges with endpoints $u$ and $v$ is deleted. We may also write, for example, that $uv\in E(G)$ when technically we should say there is an edge in $E(G)$ with endpoints $u$ and $v$.
If $\mu_G(uv)\leq 1$ for each distinct $u$ and $v$ in $V(G)$, we say that $G$ is {\em simple}.

We denote the complete graph with vertex set $V$ by $K_V$ and the complete bipartite graph with parts $U$ and $V$ by $K_{U,V}$.
If $G$ is a graph and $\l$ is a positive integer, then $\l G$ is the graph with vertex set $V(G)$ and with $\mu_{\l G} (uv)=\l \mu_{G}(uv)$ for each pair of distinct $u$ and $v$ in $V(G)$. If $H$ is a subgraph of $G$, then $G-H$ is the graph with vertex set $V(G)$ and edge set $E(G)\setminus E(H)$. Similarly, if $E\subseteq E(G)$, then $G-E$ is the graph obtained from $G$ by deleting the edges in $E$. Note that $\mu_{G-H}(uv)=\mu_G(uv)-\mu_{H}(uv)$ for each pair of distinct $u$ and $v$ in $V(H)$. Conversely, if $H$ is a graph that is edge-disjoint from $G$, then $G\cup H$ is the graph with vertex set $V(G)\cup V(H)$ and edge set $E(G)\cup E(H)$. Note that $\mu_{G\cup H}(uv)=\mu_G(uv)+\mu_{H}(uv)$ for each pair of distinct $u$ and $v$ in $V(G)\cap V(H)$.
A cycle with $m$ edges, $m\geq 2$, is called an {\em $m$-cycle} and is denoted $(v_1,v_2,\ldots,v_m)$, where $v_1,v_2,\ldots,v_m$ are the vertices of the cycle and $v_1v_2,v_2v_3,\ldots,v_{m-1}v_m,v_mv_1$ are the edges (so $2K_2$ is a $2$-cycle). A path with $m$ edges is called an {\em $m$-path} and is denoted $[v_0,v_1,\ldots,v_m]$, where $v_0,v_1,\ldots,v_m$ are the vertices of the path and $v_0v_1,v_1v_2,\ldots,v_{m-1}v_m$ are the edges. A graph is said to be {\em even} if every vertex of the graph has even degree and is said to be {\em odd} if every vertex of the graph has odd degree.

For brevity, we avoid having to make separate mention of the case where our decompositions are into cycles and a perfect matching (rather than just into cycles), as follows.  
Let $K$ be a graph and let $M=(m_1,m_2,\ldots,m_t)$ be a list of integers with $m_i\geq 2$ for $i=1,2,\ldots,t$. If $K$ is an even graph, then an $(M)$-decomposition of $K$ is a decomposition $\{G_1,G_2,\ldots,G_t\}$ such that $G_i$ is an $m_i$-cycle for $i=1,2,\ldots,t$. If $K$ is an odd graph, then an $(M)$-decomposition of $K$ is a decomposition $\{G_1,G_2,\ldots,G_t,I\}$ such that $G_i$ is an $m_i$-cycle for $i=1,2,\ldots,t$ and $I$ is a perfect matching in $K$.

A {\em packing} of a graph $K$ is a decomposition of some subgraph $G$ of $K$, and the graph $K-G$ is called the {\em leave} of the packing. An $(M)$-packing of $\l K_n$ is an $(M)$-decomposition of some subgraph $G$ of $\l K_n$ such that $G$ is an even graph if $\l(n-1)$ is even, and $G$ is an odd graph otherwise. Thus the leave of an $(M)$-packing of $\l K_n$ is an even graph and, like an $(M)$-decomposition of $\l K_n$, an $(M)$-packing of $\l K_n$ contains a perfect matching if and only if $\l (n-1)$ is odd. 

Throughout the paper we denote by $\nu_m(M)$ the number of occurrences of $m$ in the list $M$. We shall also sometimes use superscripts to specify the number of occurrences of a particular integer in a list. That is, we define $(m_1^{\alpha_1},m_2^{\alpha_2},\ldots,m_t^{\alpha_t})$ to be the list comprised of $\alpha_i$ occurrences of $m_i$ for $i=1,2,\ldots,t$. Let $M$ and $M'$ be lists of integers. It follows that for some distinct $m_1,m_2,\ldots,m_t$ we may write $M=(m_1^{\alpha_1},m_2^{\alpha_2},\ldots,m_t^{\alpha_t})$ and $M'=(m_1^{\beta_1},m_2^{\beta_2},\ldots,m_t^{\beta_t})$, where $\alpha_i,\beta_i\geq 0$ for $i=1,2,\ldots,t$. Then $\sum M=\alpha_1m_1+\alpha_2m_2+\cdots+\alpha_tm_t$, $(M,M')$ is the list $(m_1^{\alpha_1+\beta_1},m_2^{\alpha_2+\beta_2},\ldots,m_t^{\alpha_t+\beta_2})$ and, if $0\leq \beta_i\leq \alpha_i$ for $i=1,2,\ldots,t$, then $M-M'$ is the list $(m_1^{\alpha_1-\beta_1},m_2^{\alpha_2-\beta_2},\ldots,m_t^{\alpha_t-\beta_t})$.

\section{Necessity and admissible lists}\label{Section2}

For positive integers $\lambda$ and $n$, we say that a list $(m_1,m_2,\ldots,m_t)$ of integers is {\em $(\lambda,n)$-admissible} if
\begin{itemize}
\item[(A1)] $2\leq m_1,m_2,\ldots,m_t\leq n$;
\item[(A2)] $m_1+m_2+\cdots+m_t=\lambda \binom{n}{2}$ when $\lambda(n-1)$ is even;
\item[(A3)] $m_1+m_2+\cdots+m_t=\lambda \binom{n}{2}-\frac{n}{2}$ when $\lambda(n-1)$ is odd;
\item[(A4)] $\max(m_1,m_2,\ldots,m_t)+t-2\leq \frac{\lambda}{2}\binom{n}{2}$ when $\lambda$ is even; and
\item[(A5)] $\sum_{m_i=2} m_i\leq (\lambda-1)\binom{n}{2}$ when $\lambda$ is odd.
\end{itemize}

It is clear that the conditions of Theorem \ref{mainfull} are satisfied if and only if the list 
$(m_1,m_2,\ldots,m_t)$ is $(\lambda,n)$-admissible. Thus, with the above notation in hand, 
we can restate our main theorem (Theorem \ref{mainfull}) as follows.

\begin{theorem}\label{main}
For all positive integers $\lambda$ and $n$, there is an $(M)$-decomposition of $\lambda K_n$ if and only if $M$ is a $(\lambda,n)$-admissible list.
\end{theorem}

As noted above, Theorem \ref{main} is known to hold when $\l=1$ \cite{BrHoPe}, 
and we include this result here for later reference.

\begin{lemma}\label{1fold}{\rm (\cite{BrHoPe})}
Theorem \ref{main} holds for $\l=1$.
\end{lemma}

The necessity of conditions (A1)-(A3) is obvious, but proving that conditions 
(A4) and (A5) are necessary requires some work. 
The following lemma, which is a (slight) generalisation of a result in \cite{MaSm},
is used to prove the necessity of condition (A4). 

\begin{lemma}\label{BMBS}
Suppose $G$ is a  graph in which every edge has even multiplicity, $\D$ is a cycle decomposition of $G$, and $C\in \D$. Then $|\D|\leq |E(G)|/2-|E(C)|+2$.
\end{lemma}

\proof Suppose for a contradiction that there is a graph $G$ in which every edge has even multiplicity that admits a cycle decomposition $\D$ such that $|\D|> |E(G)|/2-|E(C)|+2$ for some $C\in \D$. Suppose further that, of all such graphs, $G$ has a minimal number of edges. Note that $|E(G)|>2|E(C)|$ because otherwise $|E(G)|=2|E(C)|$ and $|\D|=2$, and thus $G-C$ contains at least one edge of multiplicity at least $2$. Let $C'$ be a $2$-cycle in $G-C$, let $G'=G-C'$, let $C_1$ and $C_2$ be distinct cycles in $\D$ such that $E(C')\subset E(C_1)\cup E(C_2)$, and let $\D^*$ be a cycle decomposition of $(C_1\cup C_2)-C'$. Observe that every edge in $G'$ has even multiplicity, that $\D'=(\D\setminus\{C_1,C_2\})\cup \D^*$ is a cycle decomposition of $G'$, and that $C\in \D'$. Since $|E(G')|=|E(G)|-2$, it follows that $$|\D'|\geq (|\D|-2)+1>|E(G)|/2-|E(C)|+1=|E(G')|/2-|E(C)|+2.$$ This contradicts the minimality of $G$.\qed 

\vspace{0.3cm}

We are now ready to prove the necessity of $(\l,n)$-admissibility. 

\begin{lemma}\label{necessaryconditions}
If there is an $(M)$-decomposition of $\lambda K_n$, then $M$ is a $(\lambda,n)$-admissible list.
\end{lemma}

\proof
Let $M$ be the list $(m_1,m_2\ldots,m_t)$ and let $\D$ be an $(M)$-decomposition of $\lambda K_n$.
To show that $M$ is a $(\lambda,n)$-admissible list we need to show that conditions (A1)-(A5) hold.
It is clear that conditions (A1)-(A3) hold. 
Without loss of generality we can assume $m_1=\max(m_1,m_2,\ldots,m_t)$.
If $\lambda$ is even, then every edge in $\lambda K_n$ has even multiplicity and by Lemma \ref{BMBS}, with $G=\lambda K_n$ and $C$ a cycle of length $m_1$ in $\D$, we have $t\leq \frac{\lambda}{2}\binom{n}{2}-m_1+2$. Thus (A4) holds.
If $\lambda$ is odd, it is clear that for any pair of distinct vertices $x,y\in V(\lambda K_n)$, at most $\lambda-1$ of the edges joining $x$ and $y$ occur in $2$-cycles in $\D$. Thus (A5) holds.
\qed

Later in the paper we will often need to establish the admissibility of certain lists,
and the following lemma is a useful tool in this regard.  

\begin{lemma}\label{admissible}
Suppose $\lambda\geq 2$ and $n$ are positive integers and $M=(m_1,m_2,\ldots,m_t)$ is a list of integers satisfying
\begin{itemize}
\item[{\rm (A1)}] $2\leq m_1,m_2,\ldots,m_t\leq n$;
\item[{\rm (A2)}] $m_1+m_2+\cdots+m_t=\lambda \binom{n}{2}$ when $\lambda(n-1)$ is even; and
\item[{\rm (A3)}] $m_1+m_2+\cdots+m_t=\lambda \binom{n}{2}-\frac{n}{2}$ when $\lambda(n-1)$ is odd.
\end{itemize}
If either $\nu_2(M)< n$, or $\lambda$ is even and the two largest entries in $M$ are equal, then $M$ is $(\lambda,n)$-admissible.
\end{lemma}

\proof The result is trivially true for $n=1$, and thus we may assume that $n\geq 2$. Let $M=(m_1,m_2,\ldots,m_t)$ be a list which satisfies the conditions of the lemma. Without loss of generality we may assume that $M$ is non-increasing. By the definition of $(\lambda,n)$-admissibility, we need only show that $2\nu_2(M)\leq (\lambda-1)\binom{n}{2}$ if $\lambda$ is odd and that $m_1+t-3 < \frac{\lambda}{2}\binom{n}{2}$ if $\lambda$ is even (note that $m_1$ and $t$ are integers).\\

\noindent{\bf Case 1.} Suppose that $\lambda$ is odd. Then $\nu_2(M)<n$. Since $\lambda\geq 3$ and $n\geq 2$ we have that $2\nu_2(M)<2(n-1)\leq (\lambda-1)\binom{n}{2}$ and the result follows.\\

\noindent{\bf Case 2a.} Suppose that $\lambda$ is even and $m_1=m_2$. Then $m_1+m_2+\cdots+m_{t}=\lambda\binom{n}{2}$ and it follows that $m_3+m_4+\cdots+m_t=\lambda\binom{n}{2}-2m_1$. Thus
$$t \leq \tfrac{1}{2}\left(\lambda\textstyle{\binom{n}{2}}-2m_1\right)+2 = \tfrac{\lambda}{2}\textstyle{\binom{n}{2}}-m_1+2$$
and it follows that $m_1+t-3 < \frac{\lambda}{2}\binom{n}{2}$.\\

\noindent{\bf Case 2b.} Suppose that $\lambda$ is even and $\nu_2(M) < n$. Again $m_1+m_2+\cdots+m_{t}=\lambda\binom{n}{2}$. Let $t'=t-\nu_2(M)$. We have that  $m_2+m_3+\cdots+m_{t'}=\lambda\textstyle{\binom{n}{2}}-m_1-2\nu_2(M)$. Thus,
$$t \leq \tfrac{1}{3}\left(\lambda\textstyle{\binom{n}{2}}-m_1-2\nu_2(M)\right)+\nu_2(M)+1 = \tfrac{\lambda}{3}\textstyle{\binom{n}{2}}-\tfrac{1}{3}m_1+\tfrac{1}{3}\nu_2(M)+1.$$
So $m_1+t-3 < \frac{\lambda}{2}\binom{n}{2}$ will hold provided $\tfrac{2}{3}m_1+\tfrac{1}{3}\nu_2(M)-2 < \frac{\lambda}{6}\binom{n}{2}$ holds. Because $m_1 \leq n$, $\nu_2(M) \leq n-1$ and $\lambda \geq 2$, this latter does indeed hold. \qed

\section{Preliminary results}

In Section \ref{reduction} we will prove that it is sufficient to establish our result for a certain subset of all $(\l,n)$-admissible lists. In this section we prove a number of preliminary results which we will need for this.

\begin{lemma}\label{exactBound}
Suppose $\lambda$ and $n$ are positive integers and $M=(m_1,m_2,\ldots,m_t)$ is a $(\lambda,n)$-admissible list. If either\begin{itemize}
\item $\lambda$ is even and $\max (m_1,m_2,\ldots,m_t)+t-2=\frac{\lambda}{2}\binom{n}{2}$; or
\item $\lambda$ is odd and $\sum_{m_i=2}m_i= (\lambda-1)\binom{n}{2}$;
\end{itemize}
then there is an $(M)$-decomposition of $\lambda K_n$.
\end{lemma}

\proof Without loss of generality we may assume that $M=(m_1,m_2,\ldots,m_t)$ is non-increasing. Let $\nu=\nu_2(M)$, let $t'=t-\nu$ and let $M'=(m_1,m_2,\ldots,m_{t'})$.\\

\noindent{\bf Case 1.} Suppose that $\lambda$ is even and $m_1+t-2=\frac{\lambda}{2}\binom{n}{2}$. Since $M$ is $(\lambda,n)$-admissible, $m_1+m_2+\cdots+m_t=\lambda \binom{n}{2}$. 
It follows that
\begin{itemize}
\item $m_1+m_2+\cdots+m_{t'}=m_1+m_2+\cdots+m_{t}-2\nu=\lambda \binom{n}{2}-2\nu$ is even; and
\item $m_1+t'-2=\frac{\lambda}{2}\binom{n}{2}-\nu=\frac{1}{2}(m_1+m_2+\cdots+m_{t'})$.
\end{itemize}
It follows from Theorem $2.2$ in \cite{MaSm} that there is an $(M')$-decomposition of $2G$ for some simple graph $G$ satisfying
$$|V(G)|=\frac{1}{2}(m_1+m_2+\cdots+m_{t'})-t'+2=\tfrac{\lambda}{2}\tbinom{n}{2}-\nu-(t-\nu)+2=m_1.$$
Since $m_1\leq n$, $G$ is simple and $\lambda\geq 2$, we can relabel the vertices of $G$ so that 
$2G$ is a subgraph of $\lambda K_n$. 
It is clear that there is a $(2^\nu)$-decomposition of $\lambda K_n-2G$ and thus the result follows.\\

\noindent{\bf Case 2.} Suppose that $\lambda$ is odd and $2\nu=(\lambda-1)\binom{n}{2}$. Since $M$ is $(\lambda,n)$-admissible, $m_1+m_2+\cdots+m_t=\lambda\binom{n}{2}$ when $n$ is odd and $m_1+m_2+\cdots+m_t=\lambda\binom{n}{2}-\frac{n}{2}$ when $n$ is even. It follows that $m_1+m_2+\cdots+m_{t'}=\binom{n}{2}$ when $n$ is odd, $m_1+m_2+\cdots+m_{t'}=\binom{n}{2}-\frac{n}{2}$ when $n$ is even, and hence $M'$ is a $(1,n)$-admissible list. Thus by Lemma \ref{1fold} there is an $(M')$-decomposition of $K_n$. Furthermore, it is clear that there is a $(2^\nu)$-decomposition of $(\lambda-1)K_n$ and the result follows.\qed

The following two lemmas are taken directly from \cite{BrHoMaSm}.

\begin{lemma}\label{equalise}{\rm (\cite{BrHoMaSm})}
Let $M$ be a list of integers and let $\lambda$, $n$, $m_1$, $m_2$, $m_1'$ and $m_2'$ be positive integers such that $m_1\leq m_1'\leq m_2'\leq m_2$ and $m_1'+m_2'=m_1+m_2$. If there is an $(M,m_1,m_2)$-decomposition of $\lambda K_n$ in which an $m_1$-cycle and an $m_2$-cycle share at least two vertices, then there is an $(M,m_1',m_2')$-decomposition of $\lambda K_n$.
\end{lemma}

\begin{lemma}\label{join}{\rm (\cite{BrHoMaSm})}
Let $M$ be a list of integers and let $\lambda$, $n$, $m$, $m'$ and $h$ be positive integers such that $h\geq m+m'$ and $m+m'+h\leq n+1$. If there is an $(M,h,m,m')$-decomposition of $\lambda K_n$, then there is an $(M,h,m+m')$-decomposition of $\lambda K_n$.
\end{lemma}

In order to prove our next result we introduce the following definition.
A graph $G$ is an {\em $(a_1,a_2,\ldots,a_s)$-flower} if $G$ is the union of $s\geq 1$ cycles $A_1,A_2,\ldots,A_s$ such that
\begin{itemize}
\item $A_i$ is a cycle of length $a_i\geq 2$ for $i=1,2,\ldots,s$; and
\item if $s\geq 2$, then there is an $x\in V(G)$ such that $V(A_i)\cap V(A_j)=\{x\}$ for all $1\leq i<j\leq s$.
\end{itemize}

\begin{lemma}\label{adding3s}
Let $M$ and $M'$ be lists of integers and let $\lambda$, $n$ and $m\geq 3$  be positive integers. If there is an $(M)$-packing of $\lambda K_n$ whose leave has an $(M',m,2)$-flower as its only nontrivial component, then there is an $(M,3)$-packing of $\lambda K_n$ whose leave has an $(M',m-1)$-flower as its only nontrivial component.
\end{lemma}

\proof If $m=3$ the result is obvious. Suppose then that $m\geq 4$. Let $\P$ be an $(M)$-packing of $\lambda K_n$ which satisfies the conditions of the lemma and let $L$ be its leave. Let $A$ and $B$ be cycles in $L$ of lengths $m$ and $2$  respectively, and let $[u,v,w,x,y]$ be a path in $L$ with $u,v,w,x\in V(A)$ and $x,y\in V(B)$. By applying Lemma $2.1$ from \cite{BrHoMaSm} to $\P$ (performing the $(v,y)$-switch with origin $w$ in the terminology of that paper) we can obtain an $(M)$-packing $\P'$ of $\lambda K_n$ with a leave $L'$ such that either $L'=(L-\{vw,vu\})+\{yw,yu\}$ (if the switch has terminus $u$) or $L'=(L-\{vw,yx\})+\{yw,vx\}$ (if the switch has terminus $x$). 
In either case, it is easy to check that $\P'\cup\{(w,x,y)\}$ is an $(M,3)$-packing of $\lambda K_n$ whose leave has an $(M',m-1)$-flower as its only nontrivial component. \qed

The following lemma is a specific case of the more general Lemma 4.15 in \cite{BrHoMaSm}.

\begin{lemma}\label{addtwo3s}
Let $M$ be a list of integers and let $\lambda$ and $n\geq 5$ be positive integers. If there is an $(M,2,2,2)$-decomposition of $\lambda K_n$ in which two $2$-cycles share at least one vertex, then there is an $(M,3,3)$-decomposition of $\lambda K_n$.
\end{lemma}

We now use the above lemmas to prove some further results which will be used in Section \ref{reduction}.

\begin{lemma}\label{equalise2}
Let $M$ be a list of integers and let $\lambda$, $n$ and $m\geq 3$ be positive integers. If there is an $(M,m,2)$-decomposition of $\lambda K_n$ in which an $m$-cycle and a $2$-cycle share at least one vertex, then there is an $(M,m-1,3)$-decomposition of $\lambda K_n$.
\end{lemma}

\proof If $m=3$ the result is trivial. Suppose then that $m\geq 4$. 
If an $m$-cycle and a $2$-cycle share at least two vertices, the result follows from Lemma \ref{equalise}. 
Otherwise, there is an $(M)$-packing of $\lambda K_n$ whose leave has an $(m,2)$-floweras its only nontrivial component. Thus by Lemma \ref{adding3s} there is an $(M,3)$-packing of $\l K_n$ whose leave has an $(m-1)$-cycle as its only nontrivial component. The result follows.\qed

\begin{lemma}\label{almostall2s}
Let $M$ be a list of integers, let $\lambda$ be odd, and let $n$, $m$, $m_1'$ and $m_2'$ be positive integers satisfying $2\leq m_1'\leq m_2'\leq m$, $m_1'+m_2'=2+m$ and $2\nu_2(M)> (\lambda-1)(\binom{n}{2}-\binom{m}{2})$. If there is an $(M,m,2)$-decomposition of $\lambda K_n$, then there is an $(M,m_1,m_2)$-decomposition of $\lambda K_n$.
\end{lemma}

\proof Let $\D$ be an $(M,m,2)$-decomposition of $\lambda K_n$ an let $C$ be an $m$-cycle in $\D$. There are $\binom{n}{2}-\binom{m}{2}$ pairs of distinct vertices of $\lambda K_n$ that are not subsets of $V(C)$, and each such pair can be the vertex set of at most $(\l-1)/2$ $2$-cycles in $\D$. Since $2\nu_2(M)> (\lambda-1)(\binom{n}{2}-\binom{m}{2})$, it follows that there is a $2$-cycle in $\D$ that shares two vertices with $C$. The result then follows by Lemma \ref{equalise}.\qed

\begin{lemma}\label{two3s}
Let $M$ be a list of integers and let $\lambda$ and $n\geq 5$ be positive integers satisfying $2\nu_2(M)\geq n-5$. If there is an $(M,2,2,2)$-decomposition of $\lambda K_n$, then there is an $(M,3,3)$-decomposition of $\lambda K_n$.
\end{lemma}

\proof Let $\D$ be an $(M,2,2,2)$-decomposition of $\lambda K_n$. Since $2\nu_2(M)\geq n-5$, the number of occurrences of vertices in $2$-cycles in $\D$ is at least $n+1$ and it follows that at least two $2$-cycles share a vertex. 
The result then follows by Lemma \ref{addtwo3s}.\qed

\section{A reduction of the problem}\label{reduction}

In this section we show that to prove Theorem \ref{main} it is sufficient to prove that 
the desired decompositions exist for what we call 
$(\lambda,n)$-ancestor lists, which we now define.
For any positive integers $\lambda$ and $n$, we shall call a list $M$ a {\em $(\lambda,n)$-ancestor} list if it is $(\lambda,n)$-admissible and satisfies
\begin{itemize}
\item[(N1)] if $n=4$, then $\nu_3(M)=0$;
\item[(N2)] if $n=5$, then $\nu_3(M)+\nu_4(M)\in\{0,1\}$
\item[(N3)] if $n\geq 6$, then $\nu_4(M)+\nu_5(M)+\cdots+\nu_{n-1}(M)\in\{0,1\}$;
\item[(N4)] if $n\geq 6$ and $\nu_3(M)\geq 1$, then $\nu_{n-2}(M)+\nu_{n-1}(M)=0$; and
\item[(N5)] if $n\geq 6$ and $\nu_3(M)\geq 2$, then $\nu_2(M)\leq \lfloor\frac{n}{2}\rfloor-3$.
\end{itemize}

\begin{theorem}\label{ancestorToAll} For each pair of positive integers $\lambda$ and $n$, if there exists an $(M')$-decomposition of $\lambda K_n$ for each $(\lambda,n)$-ancestor list $M'$, then there exists an $(M)$-decomposition of $\lambda K_n$ for each $(\lambda,n)$-admissible list $M$.
\end{theorem}

\proof Let $\lambda$ and $n$ be positive integers. Throughout this proof we assume that any $(\lambda,n)$-admissible list is written in non-increasing order. For distinct $(\lambda,n)$-admissible lists $(m_1,m_2,\ldots,m_t)$ and $(m_1',m_2',\ldots,m_{t'}')$, we say the list $(m_1',m_2',\ldots,m_{t'}')$ is {\em larger} than the list  $(m_1,m_2,\ldots,m_t)$ if $t'>t$ or if $t'=t$ and $m_k'>m_k$ where $k$ is the smallest positive integer such that $m_k\ne m_k'$. Note that this defines a total order on the set of all non-increasing $(\lambda,n)$-admissible lists.

For a contradiction, suppose the theorem does not hold for $\l$ and $n$. Then there exists a largest $(\lambda,n)$-admissible list $M=(m_1,m_2,\ldots,m_t)$ such that there is no $(M)$-decomposition of $\lambda K_n$. By assumption, $M$ is not a $(\lambda, n)$-ancestor list and so at least one of the following holds.
\begin{itemize}
\item[(1)] $n=4$ and $\nu_3(M)\geq 1$.
\item[(2)] $n=5$ and $\nu_3(M)+\nu_4(M)\geq 2$.
\item[(3)] $n\geq 6$ and $\nu_4(M)+\nu_5(M)+\cdots+\nu_{n-1}(M)\geq 2$.
\item[(4)] $n\geq 6$,  $\nu_3(M)\geq 1$ and $\nu_{n-2}(M)+\nu_{n-1}(M)\geq 1$.
\item[(5)] $n\geq 6$, $\nu_3(M)\geq 2$ and $\nu_2(M)\geq \lfloor\frac{n}{2}\rfloor-2$.
\end{itemize}
Furthermore, we may assume that $\lambda\geq 2$ (by Lemma \ref{1fold}), and that $m_1+t-2<\frac{\lambda}{2}\binom{n}{2}$ if $\lambda$ is even, and $2\nu_2(M)<(\lambda-1)\binom{n}{2}$ if $\lambda$ is odd (if we have equality in either of these there exists an $(M)$-decomposition of $\lambda K_n$ by Lemma \ref{exactBound}). These strict inequalities allow us to modify the list $M$ in a number of ways to obtain a larger list which still satisfies conditions (A4) and (A5) of $(\lambda,n)$-admissibility. We now show that there exists a $(\l,n)$-admissible list $M'$ such that $M'$ is larger than $M$ and the existence of $(M')$-decomposition of $\lambda K_n$ implies the existence of an $(M)$-decomposition of $\lambda K_n$. This will suffice to complete the proof, because an $(M')$-decomposition of $\l K_n$ must exist by the maximality of $M$ and hence an $(M)$-decomposition of $\l K_n$ exists in contradiction to our assumption.

If (1) holds then $\nu_3(M)\geq 2$ (since $\sum M$ is even) and we define $M'$ to be the list obtained from $M$ by replacing two $3$'s with a $2$ and a $4$.

If (2) holds, then there exist integers $x$ and $y$ in $M$ such that $3\leq x\leq y\leq 4$, and we define $M'$ to be the list obtained from $M$ by replacing an $x$ and a $y$ with an $x-1$ and a $y+1$.

If (3) holds, then there exist integers $x$ and $y$ in $M$ such that $4\leq x\leq y\leq n-1$, and we define $M'$ as follows.
\begin{itemize}
\item[$(a)$] If $x+y\geq n+2$, then $M'$ is the list obtained from $M$ by replacing an $x$ and a $y$ with an $x-1$ and a $y+1$.
\item[$(b)$] If $x+y\leq n+1$, then $M'$ is the list obtained from $M$ by replacing an $x$ and a $y$ with an $x-2$ and a $y+2$ if  $x=4$, $\lambda$ is odd and $2\nu_2(M)= (\lambda-1)\binom{n}{2}-2$, and replacing an $x$ with a $2$ and an $x-2$ otherwise.
\end{itemize}

If (4) holds, then there exists an $x\in\{n-2,n-1\}$ in $M$ and we define $M'$ to be the list obtained from $M$ by replacing a $3$ and an $x$ with a $2$ and an $x+1$.

If (5) holds, and neither (3) nor (4) hold, we define $M'$ as follows.
\begin{itemize}
\item[$(a)$] If either $\lambda$ is even, or $\lambda$ is odd and $2\nu_2(M)\leq(\lambda-1)\binom{n}{2}-6$, then $M'$ is the list obtained from $M$ by replacing two $3$'s with three $2$'s.
\item[$(b)$] If $\lambda$ is odd and $2\nu_2(M)\geq (\lambda-1)\binom{n}{2}-4$, then $M'$ is the list obtained from $M$ by replacing two $3$'s with a $4$ and a $2$.
\end{itemize}

It is easy to see that in each case $M'$ is $(\lambda,n)$-admissible and $M'$ is larger than $M$. We now show that we can construct an $(M)$-decomposition of $\lambda K_n$ from an $(M')$-decomposition $\D$ of $\l K_n$ by applying one of Lemmas \ref{equalise}, \ref{join}, \ref{equalise2}, \ref{almostall2s} or \ref{two3s}.

If (1) holds then we can apply Lemma \ref{equalise} (since $n=4$ and thus any $4$-cycle and $2$-cycle in $\D$ must share two vertices).
If (2) holds and $(x,y)=(3,3)$, then we can apply Lemma \ref{equalise2} (since $n=5$ and thus any $4$-cycle and $2$-cycle in $\D$ must share at least one vertex). Similarly, if (2) holds and $(x,y)\ne (3,3)$, then we can apply Lemma \ref{equalise} (since $x+y\geq n+2$).
If (3) holds and $x+y\geq n+2$, then we can apply Lemma \ref{equalise}. If (3) holds and $x+y\leq n+1$, $\lambda$ is odd, $x=4$ and $2\nu_2(M)= (\lambda-1)\binom{n}{2}-2$, then we can apply Lemma \ref{almostall2s} with $m=y+2$, $m_1'=4$ and $m_2'=y$.
Otherwise, if $(3)$ holds and $x+y\leq n+1$, then we can apply Lemma \ref{join} with $m=2$, $m'=x-2$ and $h=y$.
If (4) holds, then we can apply Lemma \ref{equalise2}. If (5) holds and either $\lambda$ is even, or $\lambda$ is odd and $2\nu_2(M)\leq(\lambda-1)\binom{n}{2}-6$, then we apply Lemma \ref{two3s}.
Finally, if (5) holds, $\lambda$ is odd and $2\nu_2(M)\geq (\lambda-1)\binom{n}{2}-4$, then we can apply Lemma \ref{almostall2s} with $m=4$ and $m_1'=m_2'=3$.\qed

\section{The case $\lambda=2$}

In this section we give a proof of Theorem \ref{main} in the case $\lambda=2$. We first present two lemmas which are proved in Sections $7$ and $8$ respectively.

\begin{lemma}\label{manyHams}
If $n\geq 5$ and $M$ is a $(2,n)$-ancestor list with $\nu_n(M)>(n-3)/2$, then there is an $(M)$-decomposition of $2 K_n$.
\end{lemma}

\begin{lemma}\label{012Hams}
If $n\geq 5$ and Theorem \ref{main} holds for $2K_{n-1}$, then there is an $(M)$-decomposition of $2 K_n$ for each $(2,n)$-ancestor list $M$ satisfying $\nu_n(M)\leq (n-3)/2$.
\end{lemma}

From these lemmas we can prove the following.

\begin{lemma}\label{2fold}
Theorem \ref{main} holds for $2K_n$.
\end{lemma}

\proof
The proof is by induction on $n$. By Theorem \ref{ancestorToAll} it suffices to prove the existence of an $(M)$-decomposition of $2K_n$ for each $(2,n)$-ancestor list $M$. The result is trivial for $n\in\{1,2\}$.
If $n=3$ then $M\in\{(3,3),(2,2,2)\}$ and in each case it is clear a suitable decomposition exists.
If $n=4$ then $M\in\{(4,4,4),(4,4,2,2),(2,2,2,2,2,2)\}$ and in each case it is clear a suitable decomposition exists.
Suppose then that $n\geq 5$ and assume Theorem \ref{main} holds for all $2K_{n'}$ with $n'<n$.
Lemma \ref{manyHams} covers each $(2,n)$-ancestor list $M$ with $\nu_n(M)>(n-3)/2$, and using the inductive hypothesis, Lemma \ref{012Hams} covers those with $\nu_n(M)\leq (n-3)/2$.\qed

\section{Proof of Theorem \ref{main}}\label{SectionProof}

Lemmas \ref{1fold} and \ref{2fold} allow us to prove our main result using induction on $\lambda$. The main ingredient in the inductive step is given in the following lemma.

\begin{lemma}\label{lambdainduction}
Let $\lambda_1$, $\lambda_2$ and $n$ be positive integers such that $\lambda_1\in\{1,2\}$ and $\lambda_2$ is even. If Theorem \ref{main} holds for $\lambda_1 K_n$ and $\lambda_2 K_n$, then there is an $(M)$-decomposition of $(\lambda_1+\lambda_2) K_n$ for each $(\lambda_1+\lambda_2,n)$-ancestor list $M$ satisfying at least one of
\begin{itemize}
\item[{\rm (i)}] $\nu_n(M)\geq \lfloor\lambda_1 (\frac{n-1}{2})\rfloor$;
\item[{\rm (ii)}] $2\nu_2(M)\geq \lambda_2 \binom{n}{2}$; or
\item[{\rm (iii)}] $\nu_3(M)\geq 2$.
\end{itemize}
\end{lemma}

\proof
Let $V$ be a vertex set of size $n$, let $\sigma_1=n\lfloor\lambda_1 (\frac{n-1}{2})\rfloor$, let $\sigma_2=\lambda_2 \binom{n}{2}$, and let $M$ be a $(\lambda_1+\lambda_2,n)$-ancestor list satisfying at least one of (i), (ii) or (iii). We note that $\sum M =\sigma_1+\sigma_2$, and that for each $i\in\{1,2\}$, if $M'$ is a $(\lambda_i,n)$-admissible list, then $\sum M'=\sigma_i$.\\

\noindent{\bf Case 1.} Suppose that $M$ satisfies (i). If $\nu_2(M)<n$, let $M_1=(n^{\sigma_1/n})$ and let $M_2=M-M_1$. It follows from Lemma \ref{admissible} that $M_i$ is $(\lambda_i,n)$-admissible for each $i\in\{1,2\}$.
Thus, by assumption there is an $(M_i)$-decomposition $\D_i$ of $\lambda_i K_V$ for each $i\in\{1,2\}$. Then $\D_1\cup \D_2$ is an $(M)$-decomposition of $(\lambda_1+\lambda_2) K_V$.
If $\nu_2(M)\geq n$, let $M_1=(n^{\sigma_1/n})$ and let $M_2=(M,n,n)-(M_1,2^n)$.
It follows from Lemma \ref{admissible} that $M_i$ is $(\lambda_i,n)$-admissible for each $i\in\{1,2\}$. Thus, by assumption there is an $(M_i)$-decomposition $\D_i$ of $\lambda_i K_V$ for each $i\in\{1,2\}$. 
Let $H_1$ be an $n$-cycle in $\D_1$ and let $H_2$ be an $n$-cycle in  $\D_2$. We may assume (by relabelling vertices in $\D_2$) that there is a $(2^n)$-decomposition $\D'$ of $H_1\cup H_2$. 
Then $$(\D_1\setminus \{H_1\})\cup (\D_2\setminus \{H_2\})\cup \D'$$ is an $(M)$-decomposition of $(\lambda_1+\lambda_2) K_V$.\\

\noindent{\bf Case 2.} Suppose that $M$ satisfies (ii). Let $M_2=(2^{\sigma_2/2})$ and $M_1=M-M_2$. If $\lambda_1=1$ then $\nu_2(M)\leq \lambda_2\binom{n}{2}$ because $M$ is $(\lambda_1+\lambda_2,n)$-admissible, and thus $\nu_2(M)=\sigma_2/2$ and $\nu_2(M_1)=0$. If $\lambda_1=2$, where $M$ is the non-increasing list $(m_1,m_2,\ldots,m_t)$ say, then $m_1+t-2\leq (\frac{2+\lambda_2}{2})\binom{n}{2}$ and thus $m_1+(t-\sigma_2/2)-2\leq (\frac{2+\lambda_2}{2})\binom{n}{2}-\sigma_2/2=\binom{n}{2}$. In either case it follows, by the definition of $(\lambda,n)$-admissible and Lemma \ref{admissible}, that $M_i$ is $(\lambda_i,n)$-admissible for each $i\in\{1,2\}$.  Thus, by assumption there is an $(M_i)$-decomposition $\D_i$ of $\lambda_i K_V$ for each $i\in\{1,2\}$. Then $\D_1\cup \D_2$ is an $(M)$-decomposition of $(\lambda_1+\lambda_2) K_V$.\\

\noindent{\bf Case 3.} Suppose that $M$ satisfies (iii).
We note that if $n=3$ then $M$ also satisfies (i) and the result follows from Case $1$. Furthermore, by the properties of $(\lambda,n)$-ancestor lists, it follows that $n\notin\{1,2,4,5\}$ and thus we may assume that $n\geq 6$ and $M$ satisfies $2\nu_2(M)\leq n-6$, $\nu_4(M)+\nu_5(M)+\cdots+\nu_{n-1}(M)\in\{0,1\}$ and $\nu_{n-2}(M)+\nu_{n-1}(M)=0$.
It follows that $2\nu_2(M)+3\nu_3(M)+n\nu_n(M)\geq \sigma_1+\sigma_2-(n-3)$, and thus $3\nu_3(M)+n\nu_n(M)\geq \sigma_1+\sigma_2-(n-3)-(n-6)$. Since $\sigma_2 >2n$, we have $3\nu_3(M)+n\nu_n(M)> \sigma_1+ 9$.
If $n\nu_n(M)\geq \sigma_1$ then $M$ satisfies (i) and the result follows from Case $1$. Suppose then that $n\nu_n(M)<\sigma_1$
and hence $3\nu_3(M)> 9$.

Let $M'=M-(3,3,3)$. 
Since $n\nu_n(M)\leq \sigma_1-n$,
it follows by the definition of $(\lambda,n)$-ancestor lists that for some $\e\in \{3,4,5\}$, $M'$ can be partitioned into lists $M_1$ and $M_2$ satisfying $\sum M_1=\sigma_1-\e$ and $\sum M_2=\sigma_2-9+\e$, with $\nu_2(M_1)=0$ and $\nu_2(M_2)\leq \frac{n-6}{2}<n$. (In particular, if $M'$ is the non-increasing list $(m_1,m_2,\ldots,m_t)$ say, then a suitable partition can be obtained by packing $M_1$ with lengths $m_1$, $m_2$, $m_3$, and so on, until $\sum M_1\in\{\sigma_1-5,\sigma_1-4,\sigma_1-3\}$.) If $\e=3$, then by Lemma \ref{admissible}, $(M_1,3)$ is $(\lambda_1,n)$-admissible and $(M_2,3,3)$ is $(\lambda_2,n)$-admissible, and so the result follows by taking the union of an $(M_1,3)$-decomposition of $\lambda_1 K_V$ and an $(M_2,3,3)$-decomposition of $\lambda_2 K_V$, which exist by assumption. Suppose then that $\e\in\{4,5\}$. It follows by Lemma \ref{admissible} that $(M_1,\e)$ is $(\lambda_1,n)$-admissible and $(M_2,9-\e)$ is $(\lambda_2,n)$-admissible. Thus, by assumption there is an $(M_1,\e)$-decomposition $\D_1$ of $\lambda_1 K_V$ and an $(M_2,9-\e)$-decomposition $\D_2$ of $\lambda_2 K_V$. Let $H_1$ be an $\e$-cycle in $\D_1$ and let $H_2$ be a $(9-\e)$-cycle in $\D_2$. By relabelling vertices we may assume that $v$, $w$, $x$, $y$ and $z$ are distinct vertices in $V$ and that $\{H_1,H_2\}=\{(v,w,x,y,z),(v,w,z,x)\}$. Then
$$(\D_1\setminus \{H_1\})\cup (\D_2\setminus \{H_2\})\cup \{(v,w,x),(x,y,z),(v,w,z)\}$$
is an $(M)$-decomposition of $(\lambda_1+\lambda_2) K_V$.\qed

We now present the proof of our main Theorem.

\vspace{3mm}
\noindent{\bf Proof of Theorem \ref{main}}\\
If there is an $(M)$-decomposition of $\l K_n$, then by Lemma \ref{necessaryconditions}, $M$ is a 
$(\l,n)$-admissible list. It remains to show that if $M$ is a 
$(\l,n)$-admissible list, then there is an $(M)$-decomposition of $\l K_n$.
By Theorem \ref{ancestorToAll} we need only show there is an $(M)$-decomposition of $\lambda K_n$ for each $(\lambda,n)$-ancestor list $M$. The proof is by induction on $\lambda$. By Lemmas \ref{1fold} and \ref{2fold}, Theorem \ref{main} holds for $K_n$ and $2K_n$. So let $\lambda\geq 3$ and assume Theorem \ref{main} holds for $\lambda' K_n$ with $\lambda'<\lambda$.

Define $\lambda_1$ and $\lambda_2$ such that $\lambda_1\in\{1,2\}$ and $\lambda_2=\lambda-\lambda_1$ is even.
Let $\sigma_1=n\lfloor\lambda_1 (\frac{n-1}{2})\rfloor$ and let $\sigma_2=\lambda_2 \binom{n}{2}$. Thus $\sigma_1$ is a multiple of $n$, $\sigma_2$ is even, and $\sum M=\sigma_1+\sigma_2$. By Lemma \ref{lambdainduction}, we need only show that at least one of the following conditions holds.
\begin{itemize}
\item[(1)] $n\nu_n(M)\geq \sigma_1$.
\item[(2)] $2\nu_2(M)\geq \sigma_2$.
\item[(3)] $\nu_3(M)\geq 2$.
\end{itemize}

Note that (1) holds if $n=1$ and (2) holds if $n=2$. Similarly, if $n=3$ then $2\nu_2(M)+n\nu_n(M)=\sum M=\sigma_1+\sigma_2$ and thus (1) or (2) holds. Suppose then that $n\geq 4$. Let $k=\sum M-2\nu_2(M)-3\nu_3(M)-n\nu_n(M)$, and note by the definition of $(\lambda,n)$-ancestor lists, that $k\leq n-1$ and $k\leq n-3$ if $\nu_3(M)\geq 1$.
Suppose, for a contradiction, that none of (1), (2) or (3) holds; that is, $n\nu_n(M)<\sigma_1$, $2\nu_2(M)<\sigma_2$ and $\nu_3(M)<2$. If $\nu_3(M)=0$, then $2\nu_2(M)+3\nu_3(M)+n\nu_n(M)\leq (\sigma_2-2)+(\sigma_1-n)=\sum M-(n+2)$, and $k\geq n+2$; a contradiction. Similarly if $\nu_3(M)=1$, then $2\nu_2(M)+3\nu_3(M)+n\nu_n(M)\leq (\sigma_2-2)+3+(\sigma_1-n)=\sum M-(n-1)$, and $k\geq n-1$; a contradiction. The result follows.\qed

The remainder of the paper is devoted to filling in the details of the case $\lambda=2$ by proving Lemmas \ref{manyHams} and \ref{012Hams}.

\section{The case of more than $(n-3)/2$ Hamilton cycles}

The aim of this section is to prove Lemma \ref{manyHams} which states that, for $n\geq 5$, there is an $(M)$-decomposition of $2K_n$ for each $(2,n)$-
ancestor list $M$ satisfying $\nu_n(M) > (n-3)/2$. This will follow directly from Lemmas \ref{ManyHamsMany2s} and \ref{ManyHamsFew2s}. Throughout this section we make frequent use of circulant graphs which we define as follows. For distinct $i,j\in\{0,1,\ldots,n-1\}$, let $d_n(i,j)$ be the shortest distance from $i$ to $j$ in the $n$-cycle $(0,1,\ldots,n-1)$. If $S\subseteq \{1,2,\ldots,\lfloor n/2\rfloor\}$, then $\langle S\rangle_n$ is the simple graph with vertex set $\{0,1,\ldots,n-1\}$ and edge set  $\{\{i,j\}: d_n(i,j)\in S\}$.

\subsection{Many $2$-cycles}

In this subsection we deal with the specific case of Lemma \ref{manyHams} in which the $(2,n)$-ancestor list $M$ satisfies $\nu_2(M)\geq n/2$.

For each positive integer $n$, we define a graph $J_{n}$ by $V(J_n)=\{0,1,\ldots,n+1\}$ and $E(J_{n})=\{\{i,i+1\},\{i,i+2\}:i=0,1,\ldots,n-1\}$.
Let $M=(m_1,m_2,\ldots,m_t)$ be a list of integers with $m_i\geq 2$ for $i=1,2,\ldots,t$.
A decomposition $\{A_1,A_2,\ldots,A_t,P_1,P_2\}$ of $2J_{n}$ such that
\begin{itemize}
\item $A_i$ is a cycle of length $m_i$ for $i=1,2,\ldots,t$, and
\item $P_i$ is a path from $0$ to $n$ such that $V(P_i)=\{0,1,\ldots,n\}$ for $i=1,2$,
\end{itemize}
will be denoted $2J_{n}\rightarrow (M,n^*,n^*)$.
We note the following basic properties of $J_n$.
\begin{itemize}
\item For any integers $y$ and $n$ such that $1\leq y<n$, the graph $J_n$ is the union of $J_{n-y}$ and the graph obtained from $J_y$ by applying the vertex map $x\mapsto x+(n-y)$. Thus, if there is a decomposition $2J_{n-y}\rightarrow (M,(n-y)^*,(n-y)^*)$ and a decomposition $2J_y\rightarrow (M',y^*,y^*)$, then there is a decomposition $2J_n\rightarrow (M,M',n^*,n^*)$. We will call this construction, and the similar constructions that follow, \emph{concatenations}.
\item For $n\geq 5$, if for each $i\in\{0,1\}$ we identify vertex $i$ of $J_n$ with vertex $i+n$ of $J_n$ the resulting graph is $\langle\{1,2\}\rangle_n$. This means that for $n\geq 5$, we can obtain an $(M,n,n)$-decomposition of $2\langle\{1,2\}\rangle_n$ from a decomposition $2J_{n}\rightarrow (M,n^*,n^*)$, provided that for each $i\in\{0,1\}$ no cycle in the decomposition of $2J_n$ contains both vertex $i$ and vertex $i+n$. Note in particular that this proviso holds if the decomposition $2J_n\rightarrow (M,M',n^*,n^*)$ was formed as a non-trivial concatenation.
\end{itemize}

\begin{lemma}\label{JDecompositions1}
The following decompositions exist.
\begin{itemize}
\item[{\rm (i)}] $2J_{k}\rightarrow (k+1,2^{(k-1)/2},k^*,k^*)$, for any odd $k\geq 1$.
\item[{\rm (ii)}] $2J_{k}\rightarrow (2k_1+1,2k_2+1,2^{(k-2)/2},k^*,k^*)$, for any positive $k$, $k_1$ and $k_2$ with  $2k_1+2k_2=k$.
\end{itemize}
\end{lemma}

\proof
(i) It is easy to check that the decomposition $2J_{1}\rightarrow (2,1^*,1^*)$ exists. Let $k\geq 3$ be odd and let $A$ be the $(k+1)$-cycle on vertices $\{0,1,\ldots,k\}$ with $$E(A)=\{\{0,1\},\{k-1,k\}\}\cup\{\{i,i+2\}:i=0,1,\ldots,k-2\}.$$
Let $\P=\{A,(2,4)\}$ if $k=3$, and $\P= \{A,(k-1,k+1)\}\cup\{(i,i+1): i=2,4,\ldots,k-3\}$ otherwise. Then $\P$ is a $(k+1,2^{(k-1)/2})$-packing of $2J_k$, and in each case it is straightforward to check that the leave of $\P$ decomposes into two paths $P_1$ and $P_2$, each from $0$ to $k$, with $V(P_1)=V(P_2)=\{0,1,\ldots,k\}$ as required.

\noindent (ii) Let $k$, $k_1$ and $k_2$ be positive integers with  $2k_1+2k_2=k$. Let $A_1$ be the $(2k_1+1)$-cycle on vertices $\{0,1,\ldots,2k_1\}$ with $$E(A_1)=\{\{0,1\},\{2k_1-1,2k_1\}\}\cup\{\{i,i+2\}:i=0,1,\ldots,2k_1-2\},$$ and let $\P_1=\{A_1\}$ if $k_1=1$, and $\P_1=\{A_1\}\cup\{(i,i+1): i=2,4,\ldots,2k_1-2\}$ otherwise. 
Similarly, let $A_2$ be the $(2k_2+1)$-cycle on vertices $\{2k_1,2k_1+1,\ldots,k\}$ with $$E(A_2)=\{\{2k_1,2k_1+1\},\{k-1,k\}\}\cup\{\{i,i+2\}:i=2k_1,2k_1+2,\ldots,k-2\},$$ and let 
$\P_2=\{A_2\}$ if $k_2=1$, and $\P_2=\{A_2\}\cup\{(i,i+1): i=2k_1+1,2k_1+3,\ldots,k-3\}$ otherwise.
Then $\P=\P_1\cup \P_2 \cup \{(k-1,k+1)\}$ is a $(2k_1+1,2k_2+1,2^{(k-2)/2})$-packing of $2J_k$, and in each case it is straightforward to check that the leave of $\P$ decomposes into two paths $P_1$ and $P_2$, each from $0$ to $k$, with $V(P_1)=V(P_2)=\{0,1,\ldots,k\}$ as required.
\qed

In the following result we use concatenations of decompositions from Lemma \ref{JDecompositions1} to build decompositions of $2\langle\{1,2\}\rangle_n$.

\begin{lemma}\label{12Circ_2xn}
If $n\geq 5$ and $M=(m_1,m_2,\ldots,m_t)$ is any list satisfying
$2\leq m_i\leq n$ for $i=1,2,\ldots,t$, $2\nu_2(M)\geq n$ and $\sum M=2n$,
then there is an $(M,n,n)$-decomposition of $2\langle\{1,2\}\rangle_n$.
\end{lemma}

\proof We note that $n/2< t\leq n$ and thus $\nu_2(M)> n-t\geq 0$.
Let $M'=M-(2^{n-t})$, and let $r=t-(n-t)$ be the length of the list $M'$, noting that $r\geq 2$ and that $\sum M'=2t$. Take a partition of $M'$ in which each part is either a single even integer or a pair of odd integers. This is possible since $\sum M'$ is even and, as $M'$ contains at least one 2, this partition has at least two parts. For each part that is a single even integer, say $\{s\}$, use Lemma \ref{JDecompositions1}(i) to construct a decomposition $2J_{s-1}\rightarrow (s,2^{(s-2)/2},(s-1)^*,(s-1)^*)$, and for each part that is a pair of odd integers, say $\{s_1,s_2\}$, use Lemma \ref{JDecompositions1}(ii) to construct a decomposition
 $2J_{s_1+s_2-2}\rightarrow (s_1,s_2,2^{(s_1+s_2-4)/2},(s_1+s_2-2)^*,(s_1+s_2-2)^*)$. Concatenate all of these decompositions to obtain a decomposition
 $2J_{2t-r}\rightarrow (M',2^{t-r},(2t-r)^*,(2t-r)^*)$. Since $2t-r=n$ and $t-r=n-t$, we have $2J_n\rightarrow (M,n^*,n^*)$ and the result follows.\qed

The following lemma is from \cite{Bry1} (see Lemma 2.3 of \cite{Bry1} and its proof),
also see \cite{BryMar}.

\begin{lemma}\label{12Circ_1fold}
Let $n\geq 5$ be an integer and $M=(m_1,m_2,\ldots,m_t)$ be any list satisfying
$3\leq m_i\leq n$ for $i=1,2,\ldots,t$ and $\sum M=n$. Let $m_0=0$ and for each $j\in\{1,2,\ldots,t\}$ let $s_j=\sum_{i=0}^{j-1}m_i$ and
\begin{equation*}
C_j=
\begin{cases}
(s_j,s_j+2,\ldots,s_j+m_j-1,s_j+m_j-2,s_j+m_j-4,\ldots,s_j+1), \text{ if $m_j$ is odd;}\\
(s_j,s_j+2,\ldots,s_j+m_j-2,s_j+m_j-1,s_j+m_j-3,\ldots,s_j+1), \text{ if $m_j$ is even.}
\end{cases}
\end{equation*}
Then the leave $H$ of the packing $\mathcal{P}=\{C_1,C_2,\ldots,C_t\}$ of $\langle\{1,2\}\rangle_n$ is an $n$-cycle and so $\mathcal{P} \cup \{H\}$ is an $(M,n)$-decomposition of $\langle\{1,2\}\rangle_n$.
\end{lemma}

\begin{lemma}\label{12Circ_3xn}
If $n\geq 5$ and $M=(m_1,m_2,\ldots,m_t)$ is any list satisfying
$2\leq m_i\leq n$ for $i=1,2,\ldots,t$ and $\sum M=n$,
then there is an $(M,n,n,n)$-decomposition of $2\langle\{1,2\}\rangle_n$.
\end{lemma}

\proof
If $\nu_2(M)=0$ the result follows by combining an $(M,n)$-decomposition of $\langle\{1,2\}\rangle_n$ and an $(n,n)$-decomposition of $\langle\{1,2\}\rangle_n$, each of which exists by Lemma \ref{12Circ_1fold}.
If $\nu_2(M)=n/2$ the result follows from Lemma \ref{12Circ_2xn}. Suppose then that $1\leq \nu_2(M)<n/2$. Let $m\geq 2$ and $h\geq 0$ be integers such that $M'=M-(m,2^{2h+1})$ satisfies $\nu_2(M')=0$. (Note that $m=2$ if $\nu_2(M)$ is even, and $3 \leq m\leq n-2$ if $\nu_2(M)$ is odd.)
By Lemma \ref{12Circ_1fold} there is an $(M',m+2,4^{h},n)$-decomposition, $\D'$ say, of $\langle\{1,2\}\rangle_n$. Furthermore, we may assume that if $h\geq 1$ then $\D'$ contains the $4$-cycles $C_i=(i,i+1,i+3,i+2)$ for $i=0,4,\ldots,4(h-1)$, and that $\D'$ contains an $(m+2)$-cycle $C$ that has the path $[4h+2,4h,4h+1,4h+3]$ as a subgraph. Let $\{H_1,H_2\}$ be an $(n,n)$-decomposition of $\langle\{1,2\}\rangle_n$, with $H_1=(0,1,\ldots,n-3,n-1,n-2)$ if $n$ is even and $H_1=(0,1,\ldots,n-1)$ if $n$ is odd. In either case, $H_1$ contains the path $[0,1,\ldots,4h+3]$ (note that, if $n$ is even, $\sum M'+m\geq 4$ and hence $n \geq 4h+6$).
Then
$$\D=(\D'\setminus \{C,C_0,C_4,\ldots,C_{4(h-1)}\})\cup \{C^*,H_1^{*},H_2\}\cup \{(i,i+1):i=0,2,\ldots,4h\}$$
is an $(M,n,n,n)$-decomposition of $2\langle\{1,2\}\rangle_n$, where $C^*$ is the $m$-cycle obtained from $C$ by replacing the path $[4h+2,4h,4h+1,4h+3]$ with the path $[4h+2,4h+3]$, and $H_1^{*}$ is the $n$-cycle obtained from $H_1$ by replacing the path $[i,i+1,i+2,i+3]$ with the path $[i,i+2,i+1,i+3]$, for each $i=0,4,\ldots,4h$.
 \qed

\begin{lemma}\label{2sandns}
If $n$, $a$ and $b$ are non-negative integers satisfying $2a+2nb=n(n-5)$, then there is a $(2^a,n^{2b})$-decomposition of $2\langle \{3,4,\ldots,\lfloor n/2\rfloor\}\rangle_n$.
\end{lemma}

\proof By Lemma $3.1$ of \cite{BrMa}, there is a decomposition $\mathcal{D}$ of $\langle\{3,4,\ldots,\lfloor n/2\rfloor\}\rangle_n$ into $\lfloor\frac{n-5}{2}\rfloor$ $n$-cycles if $n$ is odd, and into $\lfloor\frac{n-5}{2}\rfloor$ $n$-cycles and a $1$-factor if $n$ is even. In either case, let $H_1,H_2,\ldots,H_b$ be distinct $n$-cycles in $\mathcal{D}$ (noting that $b \leq \lfloor\frac{n-5}{2}\rfloor$) and let $\P$ be the packing of $2\langle\{3,4,\ldots,\lfloor n/2\rfloor\}\rangle_n$ containing exactly two copies of each of $H_1,H_2,\ldots,H_b$. Clearly the leave of $\P$ is a graph in which each edge has even multiplicity. 
Thus it can be decomposed into $2$-cycles and the result follows.\qed

Finally, we have the following result.

\begin{lemma}\label{ManyHamsMany2s}
If $n\geq 5$ and $M$ is a $(2,n)$-ancestor list satisfying  $\nu_2(M)\geq n/2$ and $\nu_n(M)\geq 2$, then there is an $(M)$-decomposition of $2K_n$.
\end{lemma}

\proof
Let $M$ be a $(2,n)$-ancestor list with $\nu_2(M)\geq n/2$ and $\nu_n(M)\geq 2$. By the definition of $(2,n)$-ancestor lists it follows that
\begin{equation}\label{2andHam}
2\nu_2(M)+n\nu_n(M)\geq n(n-2).
\end{equation}
Let $b$ be the largest  integer such that $2b\leq \min(\nu_n(M)-2,n-5)$ (note that $b$ is non-negative) and let $a=n(n-5-2b)/2$. Because $2b\geq \nu_n(M)-3$ or $2b\geq n-5$, it follows from (\ref{2andHam}) that $a\leq \nu_2(M)$, and thus $(2^a,n^{2b})$ is a sublist of $M$. Now $2a+2bn=n(n-5)$ and there exists a $(2^a,n^{2b})$-decomposition of $2\langle\{3,4,\ldots,\lfloor n/2\rfloor\}\rangle_n$ by Lemma \ref{2sandns}. Note that $M'=M-(2^a,n^{2b})$ satisfies $\sum M'=4n$ and $\nu_n(M')\geq 2$. Further, since $M$ is a $(2,n)$-ancestor list, $M'$ contains either one $3$ and at most one length in $\{4,5,\ldots,n-3\}$, or no $3$'s and at most one length in $\{4,5,\ldots,n-1\}$. It follows that either $\nu_n(M')=2$, $2\nu_2(M')\geq n$ and there is an $(M')$-decomposition of $2\langle\{1,2\}\rangle_n$ by Lemma \ref{12Circ_2xn}, or $\nu_n(M')\geq 3$ and there is an $(M')$-decomposition of $2\langle\{1,2\}\rangle_n$ by Lemma \ref{12Circ_3xn}. The result follows.\qed

\subsection{Few $2$-cycles}

In this subsection we deal with the specific case of Lemma \ref{manyHams} in which the $(2,n)$-ancestor list $M$ satisfies $\nu_2(M)<n/2$.

\begin{lemma}\label{225_333} If there is an $(M,2,2,5)$-decomposition of $2K_n$ in which vertices from two vertex disjoint $2$-cycles are joined by an edge of a $5$-cycle, then there is an $(M,3,3,3)$-decomposition of $2K_n$.\end{lemma}

\proof Our aim is to show there is an $(M)$-packing of $2K_n$ whose leave admits a $(3,4,2)$-decomposition in which the $4$-cycle and $2$-cycle share at least one vertex. Then there is an $(M,3,4,2)$-decomposition of $2K_n$ in which a $4$-cycle and a $2$-cycle share at least one vertex and the result follows by Lemma \ref{equalise2}.

By our hypothesis, there is an $(M)$-packing $\P$ of $2K_n$ with leave $L$, such that $L$ admits a decomposition $\{C_1,C_2,C\}$ in which $C$ is a $5$-cycle, say $C=(u,v,x,y,z)$, and $C_1$ and $C_2$ are vertex disjoint $2$-cycles with $u\in V(C_1)$ and $v\in V(C_2)$, say $C_1=(u,u')$ and $C_2=(v,v')$. If either $y\in \{u',v'\}$ or $(u',v')=(x,z)$  then it is easy to check that the required $(3,4,2)$-decomposition of $L$ exists. Suppose then that $y\notin \{u',v'\}$ and $(u',v')\ne (x,z)$. Without loss of generality (by suitably relabelling the vertices) we may assume that $u'\notin V(C)$. By applying Lemma $2.1$ from \cite{BrHoMaSm} to $\P$ (performing the $(u',y)$-switch with origin $u$ in the terminology of that paper) we can obtain an $(M)$-packing $\P'$ of $2K_n$ with a leave $L'$ such that one of $L'=(L-\{u'u,u'u\})+\{yu,yu\}$ (if the switch has terminus $u$),  $L'=(L-\{u'u,yx\})+\{yu,u'x\}$ (if the switch has terminus $x$), or $L'=(L-\{u'u,yz\})+\{yu,u'z\}$ (if the switch has terminus $z$). 
In each case, it is easy to check that the required decomposition of $L'$ exists.\qed

\begin{lemma}\label{SumListtoMany}
Suppose $n$ and $m$ are integers with $n\geq m\geq 3$, and $M=(m_1,m_2,\ldots,m_t)$ is a list of integers satisfying $\sum M=m$ and $m_i\geq 2$ for $i=1,2,\ldots,t$.
Then there is a subgraph of $2K_n$ which admits both an $(M,n)$-decomposition and an $(m,n)$-decomposition.
\end{lemma}

\proof
If $t=1$ then $M=(m)$ and the result is obvious. Suppose then that $t\geq 2$.
Let $C_1,C_2,\ldots,C_t$ be pairwise vertex disjoint cycles, of lengths $m_1,m_2,\ldots,m_t$ respectively, in $2K_n$. For each $i=1,2,\ldots,t$, partition the edges of $C_i$ into three paths, say $P_i$, $Q_i$ and $R_i$, of lengths $1$, $m_i-2$ and $1$, respectively, and label the vertices of $C_i$ so that $P_i=[u_i,v_i]$ and $R_i=[w_i,u_i]$ (with $w_i=v_i$ if $m_i=2$). Let $L$ be the leave of the packing $\{C_1,C_2,\ldots,C_t\}$ in $2K_n$. Define paths $X_i=[u_i,v_{i+1}]$ for $i=1,2,\ldots,t-1$, $X_t=[u_t,v_1]$, $Y_i=[w_i,u_{i+1}]$ for $i=1,2,\ldots,t-2$ and $Y_{t-1}=[w_{t-1},w_t]$ in $L$ with length 1. Also define a path $Y_t$ in $L$ from $u_t$ to $u_1$ with length $n-m+1$ whose set of internal vertices is disjoint from the set $\{u_i,v_i,w_i:i=1,2,\ldots,t\}$. Observe that
\begin{itemize}
\item $C_i=P_i\cup Q_i\cup R_i$ is a cycle of length $m_i$ for each $i=1,2,\ldots,t$;
\item $H=\bigcup_{i=1}^t (Q_i\cup X_i\cup Y_i)$ is a cycle of length $n$;
\item $C=\bigcup_{i=1}^t (Q_i\cup R_i\cup X_i)$ is a cycle of length $m$; and
\item $H'=\bigcup_{i=1}^t (P_i\cup Q_i\cup Y_i)$ is a cycle of length $n$.
\end{itemize}
Then $G=\bigcup_{i=1}^t (P_i\cup 2Q_i \cup R_i \cup X_i \cup Y_i)$ is a subgraph of $2K_n$ which admits both an $(M,n)$-decomposition  $\{C_1,C_2,\ldots,C_t,H\}$, and an $(m,n)$-decomposition $\{C,H'\}$.\qed

In order to prove the main result of this subsection we require some results on decompositions of circulant graphs of the form $\langle\{n/2-1,n/2\}\rangle_{n}$, where $n$ is even. We obtain these results using graph concatenation methods similar to those in the previous subsection. Accordingly, we redefine $J_n$ to suit our purposes in this subsection.

Let $V_i=\{i,i'\}$ for each nonnegative integer $i$. Then for each even integer $n\geq 2$, we define $J_{n}$ by
$V(J_{n})=\bigcup_{i=0}^{n/2} V_i$ and $E(J_{n})=\{\{i,i+1\},\{i,i'\},\{i',(i+1)'\}:i=0,1,\ldots,n/2-1\}$.
Let $M=(m_1,m_2,\ldots,m_t)$ be a list of integers with $m_i\geq 2$ for $i=1,2,\ldots,t$.
A decomposition $\{A_1,A_2,\ldots,A_t,P_1,P_2\}$ of $J_{n}\cup I$, where $I$ is a $1$-regular graph with $V(I)=\bigcup_{i=0}^{n/2-1} V_i$, such that
\begin{itemize}
\item $A_i$ is a cycle of length $m_i$ with $V_{n/2}\cap V(A_i)=\emptyset$ for $i=1,2,\ldots,t$;
\item $P_1$ and $P_2$ are vertex disjoint paths with end vertices in $V_0\cup V_{n/2}$ such that $|E(P_1)|+|E(P_2)|=n$;
\end{itemize}
will be denoted $J_{n}\rightarrow (M,n^+)$ if each $P_i$ has one end vertex in $V_0$ and one end vertex in $V_{n/2}$, and denoted $J_{n}\rightarrow (M,n^*)$ otherwise. We note the following basic properties of $J_{n}$.
\begin{itemize}
\item For $n\geq 8$ and $n\equiv 0\pmod{4}$, if we identify vertices $0$ and $0'$ of $J_{n}$ with vertices $(n/2)'$ and $(n/2)$ respectively of $J_{n}$, the resulting graph is isomorphic to $\langle\{n/2-1,n/2\}\rangle_{n}$. Similarly, for $n\geq 6$ and $n\equiv 2\pmod{4}$, if we identify vertices
$0$ and $0'$ of $J_{n}$ with vertices $(n/2)$ and $(n/2)'$ respectively of $J_{n}$, the resulting graph is isomorphic to $\langle\{n/2-1,n/2\}\rangle_{n}$. This means that for $n\geq 6$, we can obtain an $(M,n)$-decomposition of $\langle\{n/2-1,n/2\}\rangle_{n}\cup I_i$, for some perfect matching $I_i$ in $K_{n}$, from a decomposition $J_{n}\rightarrow (M,n^*)$.
\item For any even integers $y$ and $n$ such that $2\leq y<n$, the graph $J_{n}$ is the union of $J_{n-y}$ and the graph obtained from $J_{y}$ by applying the vertex map $(x,x')\mapsto (x+(n-y)/2,(x+(n-y)/2)')$. Similarly, if $I_1$ is a $1$-regular graph with $V(I_1)=\bigcup_{i=0}^{(n-y)/2-1}V_i$ and $I_{2}$ is a $1$-regular graph with $V(I_{2})=\bigcup_{i=0}^{y/2-1}V_i$, then the union of $I_{1}$ and the graph obtained from $I_{2}$ by applying the vertex map $(x,x')\mapsto (x+(n-y)/2,(x+(n-y)/2)')$ is a $1$-regular graph $I$ with $V(I)=\bigcup_{i=0}^{n/2-1}V_i$. Thus, if there is a decomposition $J_{n-y}\rightarrow (M,(n-y)^+)$ and a decomposition $J_{y}\rightarrow (M',y^+)$, then there is a decomposition $J_{n}\rightarrow (M,M',n^+)$. Similarly, if there is a decomposition $J_{n-y}\rightarrow (M,(n-y)^+)$ and a decomposition $J_{y}\rightarrow (M',y^*)$, then there is a decomposition $J_{n}\rightarrow (M,M',n^*)$. As before, we call this method of combining decompositions \emph{concatenation}.

\end{itemize}

\begin{lemma}\label{JDecompositions}
The following decompositions exist.
\begin{itemize}
\item[{\rm (i)}] $J_{2k}\rightarrow (2k,(2k)^*)$, for any $k\geq 2$.
\item[{\rm (ii)}] $J_{2k}\rightarrow (2k_1+1,2k_2+1,(2k)^*)$, for any $k_1\geq 2$ and $k_2\geq 1$ with $k_1+k_2+1=k$.
\item[{\rm (iii)}] $J_4\rightarrow (2,2,4^*)$, $J_8\rightarrow (2,3,3,8^*)$ and $J_{12}\rightarrow (3,3,3,3,12^*)$.
\item[{\rm (iv)}] $J_{2k}\rightarrow (2k,(2k)^+)$, for any $k\geq 1$.
\item[{\rm (v)}] $J_{2k}\rightarrow (2k_1+1,2k_2+1,(2k)^+)$, for any $k_1,k_2\geq 1$ with $k_1+k_2+1=k$.
\end{itemize}
\end{lemma}

\proof In each case we give only the decomposition $\D_{2k}$ of $J_{2k}\cup I$, noting that it is then straightforward to check that (the implicitly defined) $I=(\bigcup_{G\in\D_{2k}} G)-J_{2k}$ is a $1$-regular graph with $V(I)=\bigcup_{i=0}^{k-1} V_i$ as required.\\

\noindent{(i)} Let $\D_{2k}=\{A,[0,0'],[k,k-1,\ldots,1,1',2',\ldots,k']\}$, where $A=(0,1,0',1')$ if $k=2$,

and $A=(1',0',0,1,2',2,3',3,\ldots,(k-1)',k-1)$ otherwise.\\

\noindent{(ii)} Let $\D_{2k}=\{A_1,A_2,[0,0'],P_2\}$, where
\begin{itemize}
\item $A_1=(1',0',0,1,2)$ if $k_1=2$, and\\ $A_1=(1',0',0,1,2',2,3',3,\ldots,(k_1-1)',(k_1-1),k_1)$ otherwise,
\item $A_2=((k-2)',(k-1)',k-1)$ if $k_2=1$, and\\ $A_2=(k_1',(k_1+1)',(k_1+1),(k_1+2)',(k_1+2),\ldots,(k-1)',k-1)$ otherwise,
\item $P_2=[k,k-1,\ldots,k_1,k_1',(k_1-1)',\ldots,1',1,2,\ldots,k_1-1,(k_1+1)',(k_1+2)',\ldots,k']$.
\end{itemize}

\noindent{(iii)} Let $\D_4=\{(0,1),(0',1'),[0,0'],[2,1,1',2']\}$,

let $\D_8=\{(0,1),(0',1',2'),(2,3,3'),[0,0'],[4,3,1',1,2,2',3',4']\}$, and let

$\D_{12}=\{(0,1,2),(0',1',2'),(3,3',4'),(4,5,5'),[0,0'],[6,5,1',1,3',2',2,3,4,4',5',6']\}$.\\

\noindent{(iv)} Let $\D_{2k}=\{A,[0,1,\ldots,k],[0',1',\ldots,k']\}$, where $A=(0,0')$ if $k=1$, and

 $A=(0',0,1',1,\ldots,(k-1)',k-1)$ otherwise.\\

\noindent{(v)} Let $\D_{2k}=\{A_1,A_2,P_1,P_2\}$, where
\begin{itemize}
\item $A_1= (0',0,1)$ if $k_1=1$ and \\
 $A_1=(0',0,1',1,\ldots,(k_1-1)',k_1-1,k_1)$ otherwise,
\item $A_2 = (k_1', (k-1)', k-1)$ if $k_2 = 1$ and \\
$A_2=(k_1',(k_1+1)',k_1+1,(k_1+2)',k_1+2,\ldots,(k-1)',k-1)$ otherwise,
\item $P_1 = [0,2',3',\ldots,k']$ if $k_1 = 1$ and \\
$P_1=[0,1,\ldots,k_1-1,(k_1+1)',(k_1+2)',\ldots,k']$ otherwise, and
\item $P_2=[0',1',\ldots,k_1',k_1,k_1+1,\ldots,k]$.\qed
\end{itemize}

We now use the decompositions from Lemma \ref{JDecompositions} to build larger decompositions.

\begin{lemma}\label{HamandF}
Suppose $n\geq 6$ is even and $M=(m_1,m_2,\ldots,m_t)$ is a list of integers satisfying $\sum M=n$ and $m_i\geq 2$ for $i=1,2,\ldots,t$.
Then there is a subgraph of $2K_n$ which admits both an $(M,n)$-decomposition, and a decomposition into $\langle\{n/2-1,n/2\}\rangle_n$ and a perfect matching.
\end{lemma}

\proof
Suppose first that $M=(3,3)$. Thus $n=6$ and we can easily choose $I$ such that $\langle\{2,3\}\rangle_6\cup I\cong K_6-F$ where $F$ is a perfect matching in $K_6$. The required $(3,3,6)$-decomposition then exists by Lemma \ref{1fold}. Suppose then that $M\ne (3,3)$, and note that we need only show there is a decomposition $J_n\rightarrow (M,n^*)$.
It is routine to check that, since $M$ satisfies the hypotheses of the lemma, $M$ can be written as $M=(X,Y)$ where $X$ is some (possibly empty) list and either $Y=(2k)$ for some $k \geq 2$, $Y=(2k_1+1,2k_2+1)$ for some $k_1\geq 2$ and $k_2\geq 1$, or $Y\in\{(2,2),(2,3,3),(3,3,3,3)\}$. Let $\sum Y=y$, then $J_{y}\rightarrow (Y,y^*)$ by Lemma \ref{JDecompositions} (i)--(iii). If $X$ is empty then we are finished. If $X$ is nonempty,
take a partition of $X$ in which each part is either a single even integer or a pair of odd integers. This is possible since $\sum X$ is even. For each part that is a single even integer, say $\{s\}$, use Lemma \ref{JDecompositions}(iv) to construct a decomposition $J_{s}\rightarrow (s,s^+)$, and for each part that is a pair of odd integers, say $\{s_1,s_2\}$, use Lemma \ref{JDecompositions}(v) to construct a decomposition
 $J_{s_1+s_2}\rightarrow (s_1,s_2,(s_1+s_2)^+)$. Concatenate all of these decompositions to obtain a decomposition
 $J_{n-y}\rightarrow (X,(n-y)^+)$. Then
 we can obtain the required decomposition by concatenating this decomposition with $J_{y}\rightarrow (Y,y^*)$.\qed

\begin{lemma}\label{3lists}
If $n$ is odd and $(M_1,M_2,M_3)$ is a $(2,n)$-admissible list such that
\begin{itemize}
\item $\nu_n(M_1)\geq 1$,
\item there is an $(M_1)$-decomposition of $K_n$, and
\item there is an $(M_2,\sum M_3)$-decomposition of $K_n$,
\end{itemize}
then there is an $(M_1,M_2,M_3)$-decomposition of $2K_n$.
\end{lemma}

\proof
Let $V$ be a vertex set with $|V|=n$, and let $m=\sum M_3$. Since $n\geq m\geq 3$, it follows from Lemma \ref{SumListtoMany} that there is a subgraph of $2K_V$ that admits both an $(M_3,n)$-decomposition, $\D_3$ say, and an $(m,n)$-decomposition, $\{C,H\}$ say, where $C$ is an $m$-cycle and $H$ is an $n$-cycle. Let $D_1$ be an $(M_1)$-decomposition of $K_V$ which contains the $n$-cycle $H$, and let $D_2$ be an $(M_2,m)$-decomposition of $K_V$ which contains the $m$-cycle $C$ (such decompositions can be found by taking the decompositions given by our hypotheses and relabelling vertices). Then $(\D_1\setminus\{H\})\cup (\D_2\setminus\{C\})\cup \D_3$ is an $(M_1,M_2,M_3)$-decomposition of $2K_V$ as required.\qed

\begin{lemma}\label{3listsand1factor}
If $n$ is even and $(M_1,M_2,M_3)$ is a $(2,n)$-admissible list such that
\begin{itemize}
\item there is an $(M_1)$-decomposition of $K_n-\langle \{n/2-1,n/2\}\rangle_n$,
\item there is an $(M_2)$-decomposition of $K_n$, and
\item $\nu_n(M_3)\geq 1$,
\end{itemize}
then there is an $(M_1,M_2,M_3)$-decomposition of $2K_n$.
Furthermore,
\begin{itemize}
\item[{\rm (i)}] if $\nu_2(M_3)\geq 1$ and $\nu_4(M_2)\geq 1$, then there is an $(M^\dagger)$-decomposition of $2K_n$, where $M^\dagger=(M_1,M_2,M_3,3,3)-(2,4)$; and
\item[{\rm (ii)}] if $\nu_2(M_3)\geq 2$ and $\nu_5(M_2)\geq 1$, then there is an $(M^\dagger)$-decomposition of $2K_n$, where $M^\dagger=(M_1,M_2,M_3,3,3,3)-(2,2,5)$.
\end{itemize}
\end{lemma}

\proof
Let $V$ be a vertex set with $|V|=n$ and let $\D_1$ be an $(M_1)$-decomposition of $K_V-G$, where $G\subseteq K_V$ is isomorphic to $\langle \{n/2-1,n/2\}\rangle_n$.
Let $M_3=(m_1,m_2,\ldots,m_t,n)$. Then $m_1+m_2+\ldots+m_t=n$ and by Lemma \ref{HamandF} there is an $(M_3)$-decomposition $\D_3$ of $G\cup I$, for some perfect matching $I$ in $K_V$. Observe that the cycles of lengths $m_1,m_2,\ldots,m_t$ in $\D_3$ are necessarily pairwise vertex disjoint. Finally, let $\D_2$ be an $(M_2)$-decomposition of $K_V-I$. Then $2K_V=(K_V-G)\cup (G\cup I)\cup (K_V-I)$ and $\D_1\cup \D_2\cup \D_3$ is an $(M_1,M_2,M_3)$-decomposition of $2 K_V$.

Furthermore, if $\nu_4(M_2)\geq 1$ and $\nu_2(M_3)\geq 1$, then there is a $4$-cycle $C$ in $\D_2$ and a $2$-cycle $C'$ in $\D_3$. We note that $C'$ necessarily contains an edge of $I$, say $xy$. By relabelling vertices in $\D_2$, we may assume that $x\in V(C)$. Case (i) then follows by Lemma \ref{equalise2}.

Similarly, if $\nu_5(M_2)\geq 1$ and $\nu_2(M_3)\geq 2$, then there is a $5$-cycle $C$ in $\D_2$ and two distinct $2$-cycles $C_1$ and $C_2$ in $\D_3$. We note that $C_1$ and $C_2$ contain distinct edges of $I$, say $x_1y_1$ and $x_2y_2$. By relabelling vertices in $\D_2$, we may assume that $x_1x_2\in E(C)$. Case (ii) then follows by Lemma \ref{225_333}.
\qed

Finally, we have the following result.

\begin{lemma}\label{ManyHamsFew2s}
If $n\geq 5$ and $M$ is a $(2,n)$-ancestor list satisfying $\nu_2(M)< n/2$ and $\nu_n(M)> (n-3)/2$,
then there is an $(M)$-decomposition of $2K_n$.
\end{lemma}

\proof Let $M$ be  a $(2,n)$-ancestor list with $\nu_2(M)< n/2$ and $\nu_n(M)>(n-3)/2$. By the definition of $(2,n)$-ancestor lists it follows that
\begin{equation}\label{3andHam}
3\nu_3(M)+n\nu_n(M)> n(n-3).
\end{equation}
The problem now splits according to the parity of $n$.\\

\noindent{\bf Case 1.} Suppose that $n$ is odd.
Let $M_1=(n^{(n-1)/2})$. Thus $M_1$ is a sublist of $M$ and there is an $(M_1)$-decomposition of $K_n$ by Lemma \ref{1fold}.
Furthermore, the list $M'=M-M_1$ satisfies $\sum M'=\binom{n}{2}$ and $2\nu_2(M')=2\nu_2(M)\leq (n-1)/2$. It is straightforward to show there is a sublist $M_3$ of $M'$ satisfying $\nu_2(M_3)=\nu_2(M')$ and $3\leq \sum M_3\leq n$. Let $M_2=M'-M_3$, and thus $(M_2,\sum M_3)$ is $(1,n)$-admissible and there exists an $(M_2,\sum M_3)$-decomposition of $K_n$ by Lemma \ref{1fold}. The result then follows by Lemma \ref{3lists}, noting that $M=(M_1,M_2,M_3)$ and $\nu_n(M_1)=(n-1)/2>1$.\\

\noindent{\bf Case 2.} Suppose that $n$ is even. Let $M_1=(n^{n/2-2})$. Thus $M_1$ is a sublist of $M$ and there is an $(M_1)$-decomposition of $K_n-\langle \{n/2-1,n/2\}\rangle_n$ by Lemma 3.1 of \cite{BrMa}. 
Consider the list $M'=M-M_1$. Since $2\nu_2(M')=2\nu_2(M)<n$, $n$ divides $\sum M'$ and $M'$ contains at most one length in $\{4,5,\ldots,n-1\}$, it follows that for some $\e\in\{0,1,2\}$ there is a partition of $M'$ into $M'_2$ and $M'_3$ satisfying $\sum M'_2=n(n-2)/2-\e$, $\sum M'_3=2n+\e$, $\nu_2(M'_3)=\nu_2(M')$, $\nu_3(M_3')\geq \e$ and $\nu_n(M'_3)\geq 1$. Let $M_2=(M_2',3+\e)-(3)$ and let $M_3=(M_3,2^\e)-(3^\e)$.
Thus $M_2$ is $(1,n)$-admissible and there is an $(M_2)$-decomposition of $K_n$ by Lemma \ref{1fold}. The result then follows by Lemma \ref{3listsand1factor}, noting that $M=(M_1,M_2,M_3)$ when $\e=0$, that $M=(M_1,M_2,M_3,3,3)-(2,4)$ when $\e=1$, and that $M=(M_1,M_2,M_3,3,3,3)-(2,2,5)$ when $\e=2$. \qed

\section{The case of at most $(n-3)/2$ Hamilton cycles}

The aim of this section is to prove Lemma \ref{012Hams} which states that if $n\geq 5$ and Theorem \ref{main} holds for $2K_{n-1}$, then there is an $(M)$-decomposition of $2K_n$ for each $(2,n)$-ancestor list $M$ satisfying $\nu_n(M)\leq (n-3)/2$.
We begin with the following useful lemmas.

\begin{lemma}\label{0Ham_many2s}
If there is an $(M)$-decomposition of $2K_{n-1}$, then there is an $(M,2^{n-1})$-decomposition of $2K_{n}$.
\end{lemma}

\proof
Let $U$ be a vertex set with $|U|=n-1$, let $\infty$ be a vertex not in $U$, let $V=U\cup\{\infty\}$, and let $\D$ be an $(M)$-decomposition of $2K_U$. Then $\D\cup \D_1$
is an $(M,2^{n-1})$-decomposition of $2K_V$, where
$\D_1$ is a $(2^{n-1})$-decomposition of $2K_{\{\infty\},\,U}$.\qed

\begin{lemma}\label{0Ham_few2s}
If there is an $(M,h)$-decomposition of $2K_{n-1}$, then there is an $(M,2^{n-1-h},3^h)$-decomposition of $2K_{n}$.
\end{lemma}

\proof
Let $U$ be a vertex set with $|U|=n-1$, let $\infty$ be a vertex not in $U$, let $V=U\cup\{\infty\}$ and let $\D$ be an $(M,h)$-decomposition of $2K_U$. Let $C$ be an $h$-cycle in $\D$. Then $$(\D\setminus\{C\})\cup \D_1\cup\D_2$$
is an $(M,2^{n-1-h},3^h)$-decomposition of $2K_V$, where
\begin{itemize}
\item $\D_1$ is a $(2^{n-1-h})$-decomposition of $2K_{\{\infty\},U\setminus V(C)}$; and
\item $\D_2$ is a $(3^h)$-decomposition of $2K_{\{\infty\},V(C)}\cup C$.
\end{itemize}
These decompositions are straightforward to construct.\qed

\begin{lemma}\label{cycle_flower}
Let $n\geq 5$ and $h\geq 2$ be integers. If there is an $(M)$-packing of $2K_{n}$ with leave $L$ of size $3h$ such that either
\begin{itemize}
\item[$(i)$] $L$ has a vertex of degree $2h$ and admits a decomposition into a $(2^{h-1})$-flower and an $(h+2)$-cycle, or
\item[$(ii)$] $L$ admits a decomposition into a $(2^h)$-flower and an $h$-cycle,
\end{itemize}
 then there is an $(M,3^h)$-decomposition of $2K_{n}$.
\end{lemma}

\proof
If $h=2$ and $(i)$ holds then there is an $(M,4,2)$-decomposition of $2K_n$ in which a $4$-cycle and a $2$-cycle share at least one vertex and the result follows from Lemma \ref{adding3s}.
If $h=2$ and $(ii)$ holds then there is an $(M,2,2,2)$-decomposition of $2K_n$ in which two $2$-cycles share a vertex and the result follows from Lemma \ref{addtwo3s}.  Suppose then that $h\geq 3$ and that the result holds for any $h'<h$. Our aim is to show that there is an $(M,3)$-packing of $2K_n$ with leave of size $3(h-1)$ such that either
\begin{itemize}
\item[$(a)$] $L$ has a vertex of degree $2(h-1)$ and admits a decomposition into a $(2^{h-2})$-flower and an $(h+1)$-cycle, or
\item[$(b)$] $L$ admits a decomposition into a $(2^{h-1})$-flower and an $(h-1)$-cycle.
\end{itemize}
The result will then follow by our inductive hypothesis. 
Let $V$ be a vertex set with $|V|=n$ and let $\P$ be an $(M)$-packing of $2K_V$ with leave $L$ of size $3h$ which satisfies $(i)$ or $(ii)$.\\

\noindent{\bf Case 1.} Suppose that $L$ satisfies $(i)$. Let $\{F,C\}$ be a decomposition of $L$ in which $F$ is a $(2^{h-1})$-flower and $C$ is an $(h+2)$-cycle, let $v$ be the vertex of degree $2h$ in $L$, and let $X=V(F)\cap V(C)$.
Let $[u,v,w,x,y]$ be a path in $C$ and observe that $u,w\notin X$. If $x\in X$ then $L$ contains the $3$-cycle $(v,w,x)$, and hence $\P'=\P\cup \{(v,w,x)\}$ is an $(M,3)$-packing of $2K_V$ with leave $L'$ of size $3(h-1)$. Furthermore, $L'$ decomposes into the $(2^{h-1})$-flower $F-(v,x)$ and the $(h-1)$-cycle $(C-[v,w,x])\cup[v,x]$ and thus satisfies $(b)$.
Suppose then that $x\notin X$. It follows that $|X|<h$ and hence there is a vertex $z\in V(F)\setminus X$. By applying Lemma $2.1$ from \cite{BrHoMaSm} to $\P$ (performing the $(z,x)$-switch with origin $w$ in the terminology of that paper) we can obtain an $(M)$-packing $\P'$ of $2K_V$ with leave $L'$ such that either  
$L'=L_1=(L-\{wx,yx\})+\{wz,yz\}$ (if the switch has terminus $y$), or $L'=L_2=(L-\{wx,vz\})+\{wz,vx\}$ (if the switch has terminus $v$). In either case, $\P'\cup\{(z,v,w)\}$ is an $(M,3)$-packing of $2K_V$ with leave $L'-(z,v,w)$ of size $3(h-1)$. Furthermore, $L'-(z,v,w)$ decomposes into the $(2^{h-2})$-flower $F-(z,v)$ and either the $(h+1)$-cycle $(C-[v,w,x,y])\cup [v,z,y]$ (if $L'=L_1$) or the $(h+1)$-cycle $(C-[v,w,x])\cup [v,x]$ (if $L'=L_2$), and thus satisfies $(a)$.\\

\noindent{\bf Case 2.} Suppose that $L$ satisfies $(ii)$. Let $\{F,C\}$ be a decomposition of $L$ in which $F$ is a $(2^h)$-flower and $C$ is an $h$-cycle, let $v$ be the vertex of degree $2h$ in $F$, and let $X=V(F)\cap V(C)$.
If $h=3$ then $\P\cup\{C\}$ is an $(M,3)$-packing of $2K_V$ whose leave satisfies $(b)$. Suppose then that $h\geq 4$.\\

\noindent{\bf Subcase 2a.} Suppose that $X\ne \emptyset$ and $v\notin X$. Thus there are distinct vertices $w,x,y\in V$ such that $w\in X$ and $[w,x,y]$ is a path in $C$. If $x\in X$ then $L$ contains the $3$-cycle $(v,w,x)$, and hence $\P'=\P\cup \{(v,w,x)\}$ is an $(M,3)$-packing of $2K_V$ with leave $L'$ of size $3(h-1)$. Furthermore, $L'$ decomposes into the $(2^{h-2})$-flower $F-((v,w)\cup(v,x))$ and the $(h+1)$-cycle $(C-[w,x])\cup[w,v,x]$ and thus satisfies $(a)$. Suppose then that $x\notin X$. It follows that there is a vertex $u\in V(F)\setminus (X\cup \{v\})$.  By applying Lemma $2.1$ from \cite{BrHoMaSm} to $\P$ (performing the $(u,x)$-switch with origin $w$ in the terminology of that paper) we can obtain an $(M)$-packing $\P'$ of $2K_V$ with  leave $L'$ such that either  
$L'=L_1=(L-\{wx,yx\})+\{wu,yu\}$ (if the switch has terminus $y$), or $L'=L_2=(L-\{wx,vu\})+\{wu,vx\}$ (if the switch has terminus $v$). In either case, $\P'\cup\{(u,v,w)\}$ is an $(M,3)$-packing of $2K_V$ with leave $L'-(u,v,w)$ of size $3(h-1)$. Furthermore, $L'-(u,v,w)$ decomposes into the $(2^{h-2})$-flower $F-((v,u)\cup(v,w))$ and either the $(h+1)$-cycle $(C-[w,x,y])\cup [w,v,u,y]$ (if $L'=L_1$) or the $(h+1)$-cycle $(C-[w,x])\cup [w,v,x]$ (if $L'=L_2$), and thus satisfies $(a)$.\\

\noindent{\bf Subcase 2b.} Suppose that $v\in X$. Thus there are distinct vertices $u,w,x,y\in V$ such that $u\in V(F)\setminus X$ and $[v,w,x,y]$ is a path in $C$. Note that $w\notin V(F)$. If $x\in X$ then $L$ contains the $3$-cycle $(v,w,x)$, and hence $\P'=\P\cup \{(v,w,x)\}$ is an $(M,3)$-packing of $2K_V$ with leave $L'$ of size $3(h-1)$. Furthermore, $L'$ decomposes into the $(2^{h-1})$-flower $F-(v,x)$ and the $(h-1)$-cycle $(C-[v,w,x])\cup[v,x]$ and thus satisfies $(b)$. Suppose then that $x\notin X$. By applying Lemma $2.1$ from \cite{BrHoMaSm} to $\P$ (performing the $(u,x)$-switch with origin $w$ in the terminology of that paper) we can obtain an $(M)$-packing $\P'$ of $2K_V$ with  leave $L'$ such that either  
$L'=L_1=(L-\{wx,yx\})+\{wu,yu\}$ (if the switch has terminus $y$), or $L'=L_2=(L-\{wx,vu\})+\{wu,vx\}$ (if the switch has terminus $v$). In either case, $\P'\cup\{(u,v,w)\}$ is an $(M,3)$-packing of $2K_V$ with leave $L'-(u,v,w)$ of size $3(h-1)$. Furthermore, $L'-(u,v,w)$ decomposes into the $(2^{h-1})$-flower $F-(u,v)$ and either the $(h-1)$-cycle $(C-[v,w,x,y])\cup [v,u,y]$ (if $L'=L_1$) or the $(h-1)$-cycle $(C-[v,w,x])\cup [v,x]$ (if $L'=L_2$), and thus satisfies $(b)$.\\

\noindent{\bf Subcase 2c.} Suppose that $X=\emptyset$. Thus there are distinct vertices $u,w,x,y\in V$ such that $(u,v)$ is a $2$-cycle in $F$ and $[w,x,y]$ is a path in $C$. By applying Lemma $2.1$ from \cite{BrHoMaSm} to $\P$ (performing the $(u,x)$-switch with origin $w$ in the terminology of that paper) we can obtain an $(M)$-packing $\P'$ of $2K_V$ with  leave $L'$ of size $3h$ such that either  
$L'=L_1=(L-\{wx,yx\})+\{wu,yu\}$ (if the switch has terminus $y$), or $L'=L_2=(L-\{wx,vu\})+\{wu,vx\}$ (if the switch has terminus $v$). If $L'=L_1$, then $L'$ decomposes into the $(2^{h})$-flower $F$ and the $h$-cycle $C_1=(C-[w,x,y])\cup [w,u,y]$ and, since $V(F)\cap V(C_1)=\{u\}$, we can continue as in Subcase 2a. If $L'=L_2$, then $\deg_{L'}(v)=2h$ and $L'$ decomposes into the $(2^{h-1})$-flower $F-(u,v)$ and the $(h+2)$-cycle $C_2=(C-[w,x])\cup [w,u,v,x]$ and we can continue as in Case 1.\qed

\begin{lemma}\label{fewHam_many2s}
If $n\geq 2s+3\geq 5$ and there is an $(M,(n-1)^s)$-decomposition of $2K_{n-1}$, then there is an $(M,3^s,n^s)$-packing of $2K_{n}$ whose leave has a $(2^{n-2s-1})$-flower as its only nontrivial connected component.
\end{lemma}

\proof
Let $U$ be a vertex set with $|U|=n-1$, let $\infty$ be a vertex not in $U$, let $V=U\cup\{\infty\}$, let $\D$ be an $(M,(n-1)^s)$-decomposition of $2K_U$, and let $H_1,H_2,\ldots,H_s$ be distinct $(n-1)$-cycles in $\D$.
We begin by showing there is a subgraph $G$ of $2K_U$ such that $E(G)$ contains precisely one edge from each of the cycles $H_1,H_2,\ldots,H_s$, and such that each nontrivial connected component of $G$ is a path. (As an aside, a similar result concerning the existence of such a graph $G$, in the case where $\{H_1,H_2,\ldots,H_s\}$ is a $2$-factorisation of a graph, is given in Theorem 4.5 of \cite{BlueBook}.) Construct a sequence $(G_0,U_0),(G_1,U_1),\ldots,(G_s,U_s)$, where
\begin{itemize}
\item each $G_i$ is a subgraph of $2K_U$ of size $i$ having the property that each of its nontrivial connected components is a path, and 
\item each $U_i$ is a subset of $U$ of size $i$, 
\end{itemize}
 as follows. 
Define $V(G_0)=U$, $E(G_0)=\emptyset$ and $U_0=\emptyset$. Then for each $i\in\{1,2,\ldots,s\}$ 
let $G_i$ be the graph obtained from $G_{i-1}$ by adding an edge, $x_iy_i$ say, from $E(H_i)$, such that $x_i,y_i\in U\setminus U_{i-1}$, and let $U_i$ be a subset of $U$ containing every vertex of degree 2 in $G_i$ and exactly one vertex of degree $1$ from each nontrivial connected component of $G_i$. Observe that $E(H_i)$ always contains such an edge since $V(H_i)=U$ and $|U|=n-1>2s>2|U_{i-1}|$, and that adding such an edge to $G_{i-1}$ ensures that each nontrivial connected component of $G_i$ is a path. Then $G=G_s$ is a graph with the required properties.

Let $t$ be the number of nontrivial connected components of $G$, let $p_1,p_2,\ldots,p_t$ be their respective sizes, and let $U'=U\setminus \{x_1,x_2,\ldots,x_s,y_1,y_2,\ldots,y_s\}$ (the set of vertices of degree $0$ in $G$). Observe that $t\leq s$, that $p_1+p_2+\cdots+p_t=s$, and that $|U'|=n-1-s-t\geq n-2s-1$.
 Then $$\P=(\D\setminus\{H_1,H_2,\ldots,H_s\})\cup \{H_1',H_2',\ldots,H'_s\}$$
is an $(M,n^s)$-packing of $2K_V$, where
$H_i'=(H_i-[x_i,y_i])\cup [x_i,\infty,y_i]$, for $i\in\{1,2,\ldots,s\}$.
Furthermore, the only nontrivial connected component of the leave of $\P$ is a $(p_1+2,p_2+2,\ldots,p_t+2,2^{n-1-s-t})$-flower. Let $\P_0=\P$ and for each $i=1,2,\ldots,s$, let $\P_{i}$ be the $(M,3^i,n^s)$-packing obtained by applying Lemma \ref{adding3s} to $\P_{i-1}$, choosing $m\geq 3$. Then $\P_s$ is the required packing.\qed

\begin{lemma}\label{fewHam_few2s}
If $n\geq 2s+3\geq 5$ and there is an $(M,h,(n-1)^s)$-decomposition of $2K_{n-1}$ with $h\leq n-2s-1$, then there is an $(M,2^{n-2s-1-h},3^{s+h},n^s)$-decomposition of $2K_{n}$.
\end{lemma}

\proof 
 Since $h\leq n-2s-1$, it follows from Lemma \ref{fewHam_many2s} that there is an \linebreak $(M,h,2^{n-2s-1-h},3^s,n^s)$-packing of $2K_n$ whose leave has a $(2^h)$-flower as its only nontrivial connected component. The result then follows from Lemma \ref{cycle_flower}~$(ii)$.\qed

\vspace{5mm}
\noindent{\bf Proof of Lemma \ref{012Hams}}
Let $M$ be a $(2,n)$-ancestor list with $\nu_n(M)\leq (n-3)/2$.
Since $M$ contains at most one cycle of length in $\{4,5,\ldots,n-1\}$, we have
\begin{equation}\label{2nancestorlist}
2\nu_2(M)+3\nu_3(M)+n\nu_n(M)\geq (n-1)^2.
\end{equation}

\noindent{\bf Case 1.} Suppose that $\nu_n(M)=0$. It follows from $(\ref{2nancestorlist})$ that $2\nu_2(M)+3\nu_3(M)\geq (n-1)^2$, and since $n\geq 5$, that $\nu_2(M)+\nu_3(M)\geq n$.
Let $h=0$ if $\nu_2(M)\geq n-1$, let $h=2$ if $\nu_2(M)=n-2$, let $h=n-1-\nu_2(M)$ if $\nu_2(M)\leq n-3$, and let $M'=M-(3^h,2^{n-1-h})$.
If $h=0$, then the fact that $M$ is $(2,n)$-admissible implies that $M'$ is $(2,n-1)$-admissible.
Thus, by assumption there is an $(M')$-decomposition of $2K_{n-1}$ and the result follows by Lemma \ref{0Ham_many2s}. Otherwise, $2\leq h\leq n-1$ and $\nu_2(M')\leq 1$. Then $(M',h)$ is  $(2,n-1)$-admissible (by Lemma \ref{admissible}) and by assumption there is an $(M',h)$-decomposition of $2K_{n-1}$. The result then follows by Lemma \ref{0Ham_few2s}.\\

\noindent{\bf Case 2.} Suppose that $\nu_n(M)=1$. It follows from $(\ref{2nancestorlist})$ that $2\nu_2(M)+3\nu_3(M)\geq (n-1)^2-n$, and since $n\geq 5$, that $\nu_2(M)+\nu_3(M)\geq n-1$. Furthermore, by the properties of $(2,n)$-admissible lists, it is clear that $\nu_3(M)\geq 1$.
Let $h=0$ if $\nu_2(M)\geq n-3$, let $h=2$ if $\nu_2(M)=n-4$, let $h=n-3-\nu_2(M)$ if $\nu_2(M)\leq n-5$, and let $M'=M-(n,3^{h+1},2^{n-3-h})$.
If $h=0$, then the fact that $M$ is $(2,n)$-admissible implies that $(M',n-1)$ is $(2,n-1)$-admissible. Thus, by assumption there is an $(M',n-1)$-decomposition of $2K_{n-1}$ and the result follows by Lemma \ref{fewHam_many2s} (with $s=1$). Otherwise, $2\leq h\leq n-1$ and $\nu_2(M')\leq 1$. Then $(M',h,n-1)$ is $(2,n-1)$-admissible (by Lemma \ref{admissible}) and by assumption there is an $(M',h,n-1)$-decomposition of $2K_{n-1}$. The result then follows by Lemma \ref{fewHam_few2s} (with $s=1$).\\

\noindent{\bf Case 3.} Suppose that $\nu_n(M)=s\geq 2$. If $\nu_2(M)\geq n/2$ then the required decomposition exists by Lemma \ref{ManyHamsMany2s}. Suppose then that $\nu_2(M)< n/2$. Because $s\leq (n-3)/2$, it follows from $(\ref{2nancestorlist})$ that $\nu_3(M)\geq s$. It also follows from $(\ref{2nancestorlist})$ that $2\nu_2(M)+3\nu_3(M)\geq (n-1)^2-sn$, and since $n\geq 2s+3$, that $\nu_2(M)+\nu_3(M)\geq n-s$.
Let $h=0$ if $\nu_2(M)\geq n-2s-1$, let $h=2$ if $\nu_2(M)=n-2s-2$, let $h=n-2s-1-\nu_2(M)$ if $\nu_2(M)\leq n-2s-3$, and let $M'=M-(n^s,3^{s+h},2^{n-2s-1-h})$.
If $h=0$, then $(M',(n-1)^s)$ is $(2,n-1)$-admissible (by Lemma \ref{admissible}). Thus, by assumption there is an $(M',(n-1)^s)$-decomposition of $2K_{n-1}$ and the result follows by Lemma \ref{fewHam_many2s}. Otherwise, $2\leq h\leq n-1$ and $\nu_2(M')\leq 1$. Then $(M',h,(n-1)^s)$ is $(2,n-1)$-admissible (by Lemma \ref{admissible}) and by assumption there is an $(M',h,(n-1)^s)$-decomposition of $2K_{n-1}$. The result then follows by Lemma \ref{fewHam_few2s}.\qed
\\

\noindent{\Large \bf Acknowledgements}
\vspace{3mm}

The authors acknowledge the support of the Australian Research Council via grants DP150100530, DP150100506,
DP120100790, DP120103067, DE120100040 and DP130102987.

\end{document}